\documentstyle[12pt]{article}  
\def\sq{\hbox {\rlap{$\sqcap$}$\sqcup$}}
\def\sq{\hbox {\rlap{$\sqcap$}$\sqcup$}}
\def\R{ {\rm R \kern -.31cm I \kern .15cm}}
\def\C{ {\rm C \kern -.15cm \vrule width.5pt \kern .12cm}}
\def\Z{ {\rm Z \kern -.27cm \angle \kern .02cm}}
\def\N{ {\rm N \kern -.26cm \vrule width.4pt \kern .10cm}}
\def\1{{\rm 1\mskip-4.5mu l} }
\def\lsim{\raise0.3ex\hbox{$<$\kern-0.75em\raise-1.1ex\hbox{$\sim$}}}
\def\gsim{\raise0.3ex\hbox{$>$\kern-0.75em\raise-1.1ex\hbox{$\sim$}}}
\def\noi{\noindent}

\def\beq{\begin{equation}}   \def\eeq{\end{equation}}
\def\bea{\begin{eqnarray}}  \def\eea{\end{eqnarray}}
\def\nn{\nonumber}
\def\noi{\noindent}
\def\beeq{\begin{eqnarray}} \def\eeeq{\end{eqnarray}}
\newcommand\mysection{\setcounter{equation}{0}\section}

\newcounter{hran}

\begin{document} 
\centerline{\large\bf Long range scattering for the Maxwell-Schr\"odinger system} 
 \vskip 3 truemm \centerline{\large\bf with large magnetic field data and small Schr\"odinger data\footnote{Work supported in part by NATO Collaborative Linkage Grant 979341}} 

\vskip 0.5 truecm

\centerline{\bf J. Ginibre}
\centerline{Laboratoire de Physique Th\'eorique\footnote{Unit\'e Mixte de
Recherche (CNRS) UMR 8627}}  \centerline{Universit\'e de Paris XI, B\^atiment
210, F-91405 ORSAY Cedex, France}
\vskip 3 truemm
\centerline{\bf G. Velo}
\centerline{Dipartimento di Fisica, Universit\`a di Bologna}  \centerline{and INFN, Sezione di
Bologna, Italy}

\vskip 1 truecm

\begin{abstract}
We study the theory of scattering for the Maxwell-Schr\"odinger system
in the Coulomb gauge in space dimension 3. We prove in particular
the existence of modified wave operators for that system with no size
restriction on the magnetic field data in the framework of a direct method which
requires smallness of the Schr\"odinger data, and we determine the
asymptotic behaviour in time of solutions in the range of the wave
operators.\end{abstract}

\vskip 3 truecm
\noi AMS Classification : Primary 35P25. Secondary 35B40, 35Q40, 81U99.  \par \vskip 2 truemm

\noi Key words : Long range scattering, modified wave operators, Maxwell-Schr\"odinger system.\par 
\vskip 1 truecm

\noindent LPT Orsay 04-42\par
\noindent June 2004\par

\newpage
\pagestyle{plain}
\baselineskip 18pt

\mysection{Introduction}
\hspace*{\parindent} This paper is devoted to the theory of scattering and more precisely to
the construction of modified wave operators for the
Maxwell-Schr\"odinger system (MS)$_3$ in $3 + 1$ dimensional space
time. That system describes the evolution of a charged nonrelativistic
quantum mechanical particle interacting with the (classical)
electromagnetic field it generates. It can be written as follows~: \beq
\label{1.1e} \left \{ \begin{array}{l}  i\partial_t u = - (1/2)
\Delta_A u + A_0 u \\ \\ \sq A_0  - \partial_t \left ( \partial_t A_0
+ \nabla \cdot A\right ) = |u|^2  \\ \\ \sq A + \nabla \left (
\partial_t A_0 + \nabla \cdot A\right ) = {\rm Im} \ \overline{u}
\nabla_A u \ .  \end{array} \right . \eeq

\noi Here $u$ and $(A, A_0)$ are respectively a complex valued function
and an ${I\hskip-1truemm R}^{3+1}$ valued function defined in space
time ${I\hskip-1truemm R}^{3+1}$, $\nabla_A = \nabla - iA$ , $\Delta_A =
\nabla_A^2$ and $\sq = \partial_t^2 - \Delta$ is the d'Alembertian in ${I\hskip-1truemm
R}^{3+1}$. We shall consider that system exclusively in the Coulomb
gauge $\nabla \cdot A = 0$. In that gauge, one can replace the system
(\ref{1.1e}) by a formally equivalent one in the following standard
way. The second equation  of (\ref{1.1e}) can be solved for $A_0$ by
\beq
\label{1.2e}
A_0 = - \Delta^{-1} |u|^2 = ( 4\pi |x|)^{-1} \ * \ |u|^2 \equiv g(|u|^2) 
\eeq

\noi Substituting (\ref{1.2e}) into the first and last equations of (\ref{1.1e}) yields the new system
\beq
\label{1.3e}
\left \{ \begin{array}{l}  i\partial_t u = - (1/2)
\Delta_A u + g(|u|^2) u \\ \\ \sq A  = P\ {\rm Im} \ \overline{u} \nabla_A u  \end{array} \right . 
\eeq

\noi where $P = \1 - \nabla \Delta^{-1}\nabla$ is the projector on
divergence free vector fields, together with the Coulomb gauge
condition $\nabla \cdot A = 0$ which is formally preserved by the
evolution. From now on we restrict our attention to the system
(\ref{1.3e}). \par

The (MS)$_3$ system is known to be locally well posed in sufficiently
regular spaces \cite{11r} \cite{12r} and to have global weak solutions
in the energy space \cite{9r} in various gauges including the Coulomb
gauge. However that system is so far not known to be globally well
posed in any space.\par

A large amount of work has been devoted to the theory of scattering for nonlinear equations
and systems centering on the Schr\"odinger equation, in particular for nonlinear Schr\"odinger
(NLS) equations, Hartree equations, Klein-Gordon-Schr\"odinger (KGS), Wave-Schr\"odinger (WS) 
and Maxwell-Schr\"odinger (MS) systems. As in the case of the linear Schr\"odinger
equation, one must distinguish the short range case from the long range case. In the former
case, ordinary wave operators are expected and in a number of cases proved to exist, describing
solutions where the Schr\"odinger function behaves asymptotically like a solution of the free
Schr\"odinger equation. In the latter case, ordinary wave operators do not exist and have to be
replaced by modified wave operators including a suitable phase in their definition. In that
respect, the (MS)$_3$ system (\ref{1.1e}) belongs to the borderline
(Coulomb) long range case, because of the $t^{-1}$ decay in
$L^{\infty}$ norm of solutions of the wave equation. Such is the case
also for the Hartree equation with $|x|^{-1}$ potential, for the Wave-Schr\"odinger system (WS)$_3$
in ${I\hskip-1truemm R}^{3+1}$ and for the Klein-Gordon-Schr\"odinger system (KGS)$_2$ in ${I\hskip-1truemm R}^{2+1}$.
\par

The construction of modified wave operators for the previous long range
equations and systems has been tackled by two methods. The first one
was initiated in \cite{13r} on the example of the NLS equation in
${I\hskip-1truemm R}^{1+1}$ and subsequently applied to the NLS
equation in ${I\hskip-1truemm R}^{2+1}$ and ${I\hskip-1truemm R}^{3+1}$
and to the Hartree equation \cite{1r}, to the (KGS)$_2$ system
\cite{14r} \cite{15r} \cite{16r} \cite{17r}, to the (WS)$_3$ system
\cite{18r} and to the (MS)$_3$ system \cite{19r} \cite{21r}. That
method is rather direct, starting from the original equation or system.
It will be sketched below. It is restricted to the (Coulomb) limiting
long range case, and requires a smallness condition on the asymptotic
state of the Schr\"odinger function. Early applications of the method
required in addition a support condition on the 
Fourier transform of the Schr\"odinger asymptotic state and a smallness
condition of the Klein-Gordon or Maxwell field in the case of the
(KGS)$_2$ or (MS)$_3$ system respectively \cite{14r} \cite{21r}. The
support condition was subsequently removed for the (KGS)$_2$ and
(MS)$_3$ system and the method was applied to the (WS)$_3$ system
without a support condition, at the expense of adding a correction term
to the Schr\"odinger asymptotic function \cite{15r} \cite{18r}
\cite{19r}. The smallness condition of the KG field was then removed
for the (KGS)$_2$ system, first with and then without a support
condition \cite{16r} \cite{17r}. Finally the smallness condition on the wave field was removed for the (WS)$_3$ system,
without a support condition or a correction term to the Schr\"odinger asymptotic function \cite{8r}.\par

In the present paper, we extend the results of our previous paper
\cite{8r} from the (WS)$_3$ system to the (MS)$_3$ system in the
Coulomb gauge (\ref{1.3e}). In particular we prove the existence of
modified wave operators without any smallness condition on the magnetic
potential $A$, and without a support condition or a correction term on
the asymptotic Schr\"odinger function. In addition, in the same spirit
as in \cite{8r}, we treat the problem in function spaces that are as
large as possible, namely with regularity as low as possible. As a
consequence, we require only a much lower regularity of the asymptotic
state than in previous works.\par

For completeness and although we shall not make use of that fact in the
present paper, we mention that the same problem for the Hartree
equation and for the (WS)$_3$ and (MS)$_3$ system can also be treated
by a more complex method where one first applies a phase-amplitude
separation to the Schr\"odinger function. The main interest of that
method is to remove the smallness condition on the Schr\"odinger
function, and to go beyond the Coulomb limiting case for the Hartree
equation. That method has been applied in particular to the (WS)$_3$ system and to
the (MS)$_3$ system in a special case \cite{4r} \cite{5r} \cite{6r}.
\par

We now sketch briefly the method of construction of the modified wave
operators initiated in \cite{13r}. That construction basically consists
in solving the Cauchy problem for the system (\ref{1.3e}) with infinite
initial time, namely in constructing solutions $(u,A)$ with prescribed
asymptotic behaviour at infinity in time. We restrict our attention to
time going to $+\infty$. That asymptotic behaviour is imposed in the
form of suitable approximate solutions $(u_a,A_a)$ of the system
(\ref{1.3e}). The approximate solutions are parametrized by data $(u_+,
A_+, \dot{A}_+)$ which play the role of (actually would be in simpler
e.g. short range cases) initial data at time zero for a simpler
evolution. One then looks for exact solutions $(u,A)$ of the system
(\ref{1.3e}), the difference of which with the given asymptotic ones
tends to zero at infinity in time in a suitable sense, more precisely,
in suitable norms. The wave operator is then defined traditionally as
the map $\Omega_+ : (u_+, A_+, \dot{A}_+) \to (u,A,\partial_t A)(0)$.
However what really matters is the solution $(u, A)$ in the
neighborhood of infinity in time, namely in some interval $[T, \infty
)$, and we shall restrict our attention to the construction of such
solutions. Continuing such solutions down to $t = 0$ is a somewhat
different question, connected with the global Cauchy problem at finite
times, which we shall not touch here, especially since the (MS)$_3$ system is not known to be globally well posed in any function space.\par

The construction of solutions $(u,A)$ with prescribed asymptotic
behaviour $(u_a,A_a)$ is performed in two steps. \\

\noi \underbar{Step 1}. One looks for $(u,A)$ in the form $(u,A) = (u_a+v, A_a + B)$ with $\nabla \cdot A_a = \nabla \cdot B = 0$. The system satisfied by the new functions $(v,B)$ can be written as
\beq \label{1.4e}
\left \{ \begin{array}{l}  i\partial_t v = - (1/2)
\Delta_A v + g(|u|^2)v + G_1 - R_1 \\ \\ \sq B =  G_2 - R_2   \end{array} \right . 
 \eeq

\noi where $G_1$ and $G_2$ are defined by
\beq \label{1.5e}
\left \{ \begin{array}{l}  G_1= iB \cdot \nabla_{A_a} u_a + (1/2)
B^2 u_a + g\left ( |v|^2 + 2\ {\rm Re}\ \overline{u}_a v \right ) u_a \\ \\ G_2 = P\ {\rm Im}\left ( \overline{v} \nabla_A v + 2 \overline{v}\nabla_A u_a\right )- P\ B |u_a|^2  \end{array} \right . 
 \eeq

\noi and the remainders are defined by
\beq \label{1.6e}
\left \{ \begin{array}{l} R_1 =  i\partial_t u_a + (1/2)
\Delta_{A_a} u_a - g \left ( |u_a|^2\right ) u_a \\ \\ R_2 = \sq A_a  - P \ {\rm Im} \ \overline{u}_a \nabla_{A_a} u_a \ . \end{array} \right . 
 \eeq
 
 \noi It is technically useful to consider also the partly linearized system for functions $(v', B')$
\beq \label{1.7e}
\left \{ \begin{array}{l}  i\partial_t v' = - (1/2)
\Delta_A v' + g(|u|^2)v' + G_1 - R_1 \\ \\ \sq B' =  G_2 - R_2   \ .\end{array} \right . 
 \eeq
 
\noi The first step of the method consists in solving the system
(\ref{1.4e}) for $(v, B)$, with $(v, B)$ tending to zero at infinity in
time in suitable norms, under assumptions on $(u_a, A_a)$ of a general
nature, the most important of which being decay assumptions on the
remainders $R_1$ and $R_2$. That can be done as follows. One first
solves the linearized system (\ref{1.7e}) for $(v',B')$ with given $(v,
B)$ and initial data $(v', B')(t_0) = 0$ for some large finite $t_0$.
One then takes the limit $t_0 \to \infty$ of that solution, thereby
obtaining a solution $(v',B')$ of (\ref{1.7e}) which tends to zero at
infinity in time. That construction defines a map $\phi : (v,B) \to
(v',B')$. One then shows by a contraction method that the map $\phi$
has a fixed point. That first step will be performed in Section 2. \\

\noi \underbar{Step 2.} The second step of the method consists in
constructing approximate asymptotic solutions $(u_a, A_a)$ satisfying
the general estimates needed to perform Step~1. With the weak time
decay allowed by our treatment of Step 1, one can take the simplest
version of the asymptotic form used in previous works \cite{6r}
\cite{19r} \cite{21r}. Thus we choose 
\beq
\label{1.8e}
u_a = MD \exp (- i \varphi ) w_+
\eeq 

\noi where
\beq
\label{1.9e}
M \equiv M(t) = \exp (ix^2/2t) \ ,
\eeq
\beq
\label{1.10e} 
D(t) = (it)^{-n/2} D_0(t) \quad , \quad ( D_0 (t) f) (x) = f(x/t) \ ,
\eeq

\noi $\varphi$ is a real phase to be chosen below and $w_+ = Fu_+$. We
furthermore choose $A_a$ in the form $A_a = A_0 + A_1$
where $A_0$ is the solution of the free wave equation $\sq A_0 = 0$
given by 
\beq
\label{1.11e}
A_0 = \cos \omega t \ A_+ + \omega^{-1} \sin \omega t \ \dot{A}_+
\eeq

\noi where $\omega = (- \Delta)^{1/2}$, and where 
\beq
\label{1.12e}
A_1(t) = \int_t^{\infty} dt' (\omega t')^{-1} \sin (\omega (t'-t))P\ x |u_a (t')|^2 \ .
\eeq

\noi Substituting (\ref{1.8e}) into (\ref{1.12e}) yields
\beq
\label{1.13e}
A_1(t) = t^{-1} D_0(t) \widetilde{A}_1
\eeq

\noi where
\beq
\label{1.14e}
\widetilde{A}_1 =  \int_1^{\infty} d\nu \ \nu^{-3} \omega^{-1} \sin (\omega (\nu - 1)) D_0 (\nu ) P\ x|w_+|^2 \ .
\eeq

\noi In particular $\widetilde{A}_1$ is constant in time. We finally
choose $\varphi$ by imposing 
\beq
\label{1.15e}
\varphi (1) = 0\quad , \qquad \partial_t \varphi = t^{-1}\left ( g(|w_+|^2) - x \cdot \widetilde{A}_1 \right )
\eeq

\noi so that 
\beq
\label{1.16e}
\varphi = (\ell n\ t) \left ( g(|w_+|^2) - x\cdot \widetilde{A}_1 \right )\ .
\eeq

We shall show in Section 3 that the previous choice fulfills the
conditions needed for Step 1, under suitable assumptions on the
asymptotic state $(u_+, A_+, \dot{A}_+)$.\par

In order to state our results we introduce some notation. We denote by
$F$ the Fourier transform, by $<\cdot , \cdot >$ the scalar product in $L^2$ and by $\parallel \cdot \parallel_r$ the norm
in $L^r \equiv L^r ({I\hskip-1truemm R}^3)$, $1 \leq r \leq \infty$ and we define $\delta (r) = 3/2 - 3/r$.
For any nonnegative integer $k$ and for $1 \leq r \leq \infty$, we
denote by $W_r^k$ the Sobolev spaces 
$$W_r^k = \left \{ u : \parallel u; W_r^k\parallel \ = \sum_{\alpha : 0 \leq |\alpha | \leq k} \parallel \partial_x^{\alpha} u \parallel_r \ < \infty \right \}$$

\noi where $\alpha$ is a multiindex, so that $H^k = W_2^k$. We shall
need the weighted Sobolev spaces $H^{k,s}$ defined for $k$, $s \in
{I\hskip-1truemm R}$ by 
$$H^{k,s} = \left \{ u : \parallel u; H^{k,s}\parallel \ = \ \parallel (1 + x^2)^{s/2} (1 - \Delta )^{k/2} u \parallel_2 \ < \infty \right \}$$

\noi so that $H^k = H^{k,0}$. For any interval $I$, for any Banach
space $X$ and for any $q$, $1 \leq q \leq \infty$, we denote by $L^q(I,
X)$ (resp. $L_{loc}^q(I,X)$) the space of $L^q$ integrable (resp.
locally $L^q$ integrable) functions from $I$ to $X$ if $q < \infty$ and
the space of measurable essentially bounded (resp. locally essentially
bounded) functions from $I$ to $X$ if $q = \infty$. For any $h \in
{\cal C} ([1, \infty ), {I\hskip-1truemm R}^+)$, non increasing and
tending to zero at infinity and for any interval $I \subset [1, \infty
)$, we define the space 

$$X(I) = \Big \{ (v, B):v\in {\cal C}(I, H^2) \cap {\cal C}^1(I,L^2),$$
$$\parallel (v, B);X(I)\parallel\ \equiv \ \mathrel{\mathop {\rm Sup}_{t \in I }}\ h(t)^{-1} \Big ( \parallel v(t);H^2\parallel \ + \ \parallel \partial_t v(t) \parallel_2 \ + \ \parallel v; L^{8/3}(J,W_4^1)\parallel$$
\beq
\label{1.17e}
+ \ \parallel B;L^4(J,W_4^1)\parallel \ + \ \parallel \partial_t B;L^4 (J, L^4)\parallel \Big ) < \infty \Big \} \eeq

\noi where $J = [t, \infty) \cap I$. \par

We can now state our result. \\

\noi{\bf Proposition 1.1.} {\it Let $h(t) = t^{-1}(2 + \ell n \ t
)^2$ and let $X(\cdot )$ be defined by (\ref{1.17e}). Let $u_a$ be
defined by (\ref{1.8e}) with $w_+ = Fu_+$ and with $\varphi$ defined
by (\ref{1.16e}) (\ref{1.2e}) (\ref{1.14e}). Let $A_a = A_0 + A_1$
with $A_0$ defined by (\ref{1.11e}) and $A_1$ by (\ref{1.13e})
(\ref{1.14e}). Let $u_+ \in H^{3,1}\cap H^{1,3}$ with $\parallel xw_+
\parallel_4$ and $\parallel w_+ \parallel_3$ sufficiently small. Let
$\nabla^2 A_+$, $\nabla \dot{A}_+$, $\nabla^2(x\cdot A_+)$ and $\nabla
(x \cdot \dot{A}_+) \in W_1^1$ with $A_+$, $x \cdot A_+ \in L^3$ and
$\dot{A}_+$, $x \cdot \dot{A}_+ \in L^{3/2}$ and let $\nabla \cdot A_+
= \nabla \cdot \dot{A}_+ = 0$. \par

Then there exists $T$,
$1 \leq T < \infty$ and there exists a unique solution $(u,A)$ of the
system (\ref{1.3e}) such that $(u - u_a, A-A_a) \in X([T, \infty
))$. Furthermore $\nabla (A - A_a)$, $\partial_t (A - A_a) \in {\cal C}([T, \infty ), L^2)$ and $A$ satisfies the estimate
\beq
\label{1.18e}
\parallel \nabla (A - A_a)(t) \parallel_2\ \vee \  \parallel \partial_t (A - A_a)(t) \parallel_2\ \leq C\ t^{-3/2} (2 + \ell n\ t)^2
\eeq

\noi for some constant $C$ depending on $(u_+, A_+, \dot{A}_+)$ and for all $t \geq T$.}\\

\noi {\bf Remark 1.1} The only smallness conditions bear on $\parallel x w_+ \parallel_4$ and
on $\parallel w_+ \parallel_3$ and are required by the magnetic interaction and the Hartree interaction (\ref{1.2e}) respectively. 
In particular there is no smallness condition on $(A_+, \dot{A}_+)$.\\

\noi {\bf Remark 1.2.} The assumptions $A_+$, $x\cdot A_+ \in L^3$ and
$\dot{A}_+$, $x \cdot \dot{A}_+ \in L^{3/2}$ serve to exclude the
occurence of constant terms in $A_+$, $x \cdot A_+$, $\dot{A}_+$, $x
\cdot \dot{A}_+$ and of terms linear in $x$ in $A_+$, $x \cdot A_+$,
but are otherwise implied by the $W_1^1$ assumptions on those
quantities through Sobolev inequalities.\\

\noi {\bf Remark 1.3.} The assumptions on $A_+$, $\dot{A}_+$ imply that
$\omega^{1/2}A_+$, $\omega^{-1/2} \dot{A}_+ \in H^1$ through Sobolev
inequalities, and therefore also that $\nabla A_+$, $\dot{A}_+ \in
L^2$. As a consequence the free wave solution $A_0$ defined by
(\ref{1.11e}) belongs to $L^4({I\hskip-1truemm R}, W_4^1)$ by Strichartz inequalities,
with $\partial_t A_0 \in L^4({I\hskip-1truemm R},L^4)$ \cite{3r}. In particular $A_0$
satisfies the local in time regularity of $B$ required in the
definition of the space $X(\cdot)$. Furthermore $\nabla A_0$,
$\partial_t A_0 \in ({\cal C} \cap L^{\infty})({I\hskip-1truemm R},L^2)$, namely $A_0$ is
a finite energy solution of the wave equation.

\mysection{The Cauchy problem at infinite initial time}
\hspace*{\parindent} In this section we perform the first step of the
construction of solutions of the system (\ref{1.3e}) as described in
the introduction, namely we construct solutions $(v, B)$ of the system
(\ref{1.4e}) defined in a neighborhood of infinity in time and tending
to zero at infinity under suitable regularity and decay assumptions on
the asymptotic functions $(u_a, A_a)$ and on the remainders $R_i$. 
As a preliminary to that study, we need to solve the Cauchy problem
with finite initial time for the linearized system (\ref{1.7e}). That
system consists of two independent equations. The second one is simply
a wave equation with an inhomogeneous term and the Cauchy problem with
finite or infinite initial time for it is readily solved under suitable
assumptions on the inhomogeneous term, which will be fulfilled in the
applications. The first one is a Schr\"odinger equation with time
dependent magnetic and scalar potentials and with time dependent inhomogeneity, which we
rewrite in a more concise form and with slightly different notation as
\beq \label{2.1e} i\partial_t v = - (1/2) \Delta_A v + Vv + f \ . 
\eeq

We first give some preliminary results on the Cauchy problem with
finite initial time for that equation at the level of regularity of
$H^2$. The following proposition is a minor variation of Proposition 3.2 in \cite{7r}.\\

\noi {\bf Proposition 2.1.} {\it Let $I$ be an interval, let $A \in {\cal C}(I, L^4 +
L^{\infty})$, $\partial_t A \in L_{loc}^1(I, L^4 + L^{\infty})$, $V \in
{\cal C}(I,L^2 + L^{\infty})$, $\partial_t V \in L_{loc}^1(I, L^2 +
L^{\infty})$, $f \in {\cal C}(I,L^2)$ and $\partial_t f \in L^1_{loc}(I,
L^2)$. Let $t_0 \in I$ and $v_0 \in H^2$. Then\par

(1) There exists a unique solution $v \in {\cal C}(I, H^2) \cap {\cal C}^1(I,L^2)$ of (\ref{2.1e}) in $I$ with $v(t_0) = v_0$.
That solution is actually unique in ${\cal C}(I,H^1)$. For all $t \in I$, the following equality holds~: 
\beq \label{2.2e} 
\parallel v(t)\parallel_2^2\ -
\ \parallel  v_0\parallel_2^2\ = \int_{t_0}^t dt'\ 2\ {\rm Im}\
<v,f> (t')\ . 
\eeq

(2) Let in addition $A \in L_{loc}^2(I, L^{\infty})$,
$\nabla A \in L_{loc}^1(I, L^{\infty})$ and $V \in L_{loc}^1(I, L^{\infty})$. Then for all
$t\in I$, the following equality holds~: 
\beq \label{2.3e} 
\parallel \partial_t v(t)\parallel_2^2\ - \ \parallel  \left ( - (1/2) \Delta_A v_0 + V v_0 + f\right ) (t_0)\parallel_2^2\ =
\int_{t_0}^t dt'\ 2\ {\rm Im}\ <\partial_t v,f_1>(t') 
\eeq

\noi where
\beq \label{2.4e} 
f_1 = i \left ( \partial_t A\right )\cdot \nabla_Av + \left ( \partial_t V\right ) v + \partial_t f \ .
\eeq

\noi Furthermore the solution is unique in ${\cal C}(I,L^2)$. }\\

We shall make an essential use of the well-known Strichartz
inequalities for the Schr\"odinger equation \cite{2r} \cite{10r}
\cite{22r}, which we recall for completeness. We define
\beq \label{2.5e} U(t) = \exp (i(t/2) \Delta ) \ . 
\eeq

\noi A pair of H\"older exponents $(q, r)$ will be called admissible if
$0 \leq 2/q\ = 3/2 - 3/r \leq 1$. For any $r$, $1 \leq r \leq \infty$, we
define $\overline{r}$ by $1/r + 1/\overline{r} = 1$.\\

\noi {\bf Lemma 2.1.} {\it The following inequalities hold. \par

(1) For any admissible pair $(q, r)$ and for any $u \in L^2$ 
\beq
\label{2.6e} \parallel U(t) u ; L^q ({I\hskip-1truemm R}, L^r)
\parallel \ \leq \ C \parallel u \parallel_2 \ . 
\eeq

(2) Let $I$ be an interval and let $t_0 \in I$. Then for any admissible
pairs $(q_i, r_i)$, $i = 1,2$,} 
\beq \label{2.7e} \parallel 
\int_{t_0}^t dt' \ U(\cdot - t') f(t') ; L^{q_1}(I,L^{r_1})\parallel \
\leq \ C \parallel f;L^{\overline{q}_2}(I, L^{\overline{r}_2})\parallel
\ . 
\eeq

In addition to the Strichartz inequalities for the Schr\"odinger
equation, we shall need special cases of the Strichartz
inequalities for the wave equation \cite{3r} \cite{10r}. Let $I$ be an
interval, let $t_0 \in I$ and let $B(t_0) = \partial_t B(t_0) = 0$. Then
\beq
\label{2.8e}
\parallel B;L^4(I,L^4) \parallel \ \leq \ C \parallel  \sq B; L^{4/3} (I, L^{4/3})\parallel  \ ,
\eeq
\beq
\label{2.9e}
\parallel \nabla B;L^4(I,L^4) \parallel \ \vee \ \parallel \partial_t B ;L^4(I,L^4) \parallel \ \leq \ C \parallel  \nabla \sq B; L^{4/3} (I, L^{4/3})\parallel  \ ,
\eeq
\beq
\label{2.10e}
\mathrel{\mathop {\rm Sup}_{t\in I}}\ \left ( \parallel \nabla B(t) \parallel_2 \ \vee \ \parallel \partial_t B(t) \parallel_2 \right ) \ \leq \ \parallel \sq B; L^1(I,L^2)\parallel \ .
\eeq

We now begin the construction of solutions of the system (\ref{1.4e}).
For any $T$, $t_0$ with $1 \leq T < t_0 \leq \infty$, we denote by $I$
the interval $I = [T, t_0]$ and for any $t \in I$, we denote by $J$ the
interval $J = [t, t_0]$. In all this section, we denote by $h$ a
function in ${\cal C}([1, \infty ), {I\hskip-1truemm R}^+)$ such that
for some $\lambda > 0$, the function $\overline{h}(t) \equiv
t^{\lambda}h(t)$ is non increasing and tends to zero as $t \to \infty$, and we denote by $j$, $k$ nonnegative integers.
\par

We shall make repeated use of the following lemma. \\

\noi {\bf Lemma 2.2.} {\it Let $1 \leq q$, $q_k \leq \infty$ ($1 \leq k \leq n$) be such that 
$$\mu \equiv 1/q - \sum_k 1/q_k \geq 0 \ .$$

\noi Let $f_k \in L^{q_k}(I)$ satisfy
\beq
\label{2.11e}
\parallel f_k;L^{q_k}(J) \parallel\ \leq N_k\ h(t)
\eeq

\noi for $1 \leq k \leq n$, for some constants $N_k$ and for all $t \in I$. \par

Let $\rho \geq 0$ such that $n\lambda + \rho > \mu$. Then the following inequality holds for all $t\in I$}
\beq
\label{2.12e}
\parallel \left ( \prod_k f_k \right ) t^{-\rho} ; L^q(J) \parallel\ \leq C \left ( \prod_k N_k\right ) h(t)^n \ t^{\mu - \rho}
\eeq

\noi where
\beq
\label{2.13e}
C = \left ( 1 - 2^{-q(n\lambda + \rho - \mu )}\right )^{-1/q} \ .
\eeq
\vskip 5 truemm

\noi {\bf Proof.} For $t \in I$, we define $I_j = [t2^j, t2^{j+1}] \cap
I$ so that $J = \ \displaystyle{\mathrel{\mathop {\cup}_{j\geq 0}}}\
I_j$. We then rewrite $L^q(J) = \ell_j^q(L^q(I_j))$. We estimate
$$\parallel \left ( \prod_k f_k \right ) t^{-\rho} ; L^q(J) \parallel\ \leq\ \parallel \left ( \prod_k \parallel f_k ;L^{q_k}(I_j)\parallel\right ) \parallel t^{-\rho} ; L^{1/\mu}(I_j) \parallel \ ; \ell_j^q \parallel$$
$$\leq \ \left ( \prod_k N_k \right ) \parallel h(t2^j)^n (t2^j)^{-\rho + \mu}; \ell_j^q \parallel$$
$$\leq \ \left ( \prod_k N_k \right ) \overline{h}(t)^n \ t^{-n\lambda - \rho + \mu} \parallel 2^{j(-n\lambda - \rho + \mu )};\ell_j^q \parallel$$
 
\noi from which (\ref{2.12e}) follows. \par \nobreak \hfill $\sq$ \\

\noi {\bf Remark 2.1.} In some special cases, the dyadic decomposition
is not needed for the proof of Lemma 2.1. For instance if all the $q_k$
are infinite, one can estimate
\bea
\label{2.14e}
&&\parallel h(t)^n \ t^{-\rho} \parallel_q \ \leq \overline{h}(t)^n \parallel t^{-\rho - n \lambda} \parallel_q\nn \\
&&\leq C\ \overline{h}(t)^n\ t^{-\rho - n \lambda + 1/q} = C\ h(t)^n\ t^{-\rho + \mu}
\eea

\noi by a direct application of H\"older's inequality in $J$. The same situation occurs if $\rho > \mu$.\\

In order to estimate the Hartree interaction term (\ref{1.2e}), we shall use the following Lemma. We recall that $\delta (r) = 3/2 - 3/r$.\\

\noi {\bf Lemma 2.3.} {\it The following estimates hold. \par
(1) 
\beq
\label{2.15e}
\parallel g(\overline{v}_1 v_2)v_3 \parallel_{\overline{r}_4} \ \leq \ C \prod_{1\leq i \leq 3} \parallel v_i \parallel_{r_i}
\eeq

\noi for $0 \leq \delta_i = \delta (r_i) \leq 1$, $1 \leq i \leq 4$, $\sum \delta_i = 1$, $0 < \delta_1 + \delta_2 < 1$.\par

(2)
\beq
\label{2.16e}
\parallel g(\overline{v}_1 v_2) \parallel_{\infty} \ \leq \ C  \parallel v_2 \parallel_{r_2} \left ( \parallel v_1 \parallel_{r_{1+}} \ \parallel v_1 \parallel_{r_{1-}}\right )^{1/2} 
\eeq

\noi for $0 < 3/r_1 = 2 - 3/r_2 \leq 2$, $1/r_{1\pm} = (1 \mp \varepsilon )/r_1$, $\varepsilon > 0$.}\\

\noi {\bf Proof.} \underline{Part (1)} follows from the H\"older and Hardy-Littlewood-Sobolev inequalities. \\

\noi \underline{Part (2)} is proved by separating $|x|^{-1}$ into short
and long distance parts, applying the H\"older inequality, and
optimizing the result with respect to the point of separation (see
\cite{1r}). \par \nobreak \hfill $\sq$ \par

We can now state the main result of this section.\\

\noi {\bf Proposition 2.2.} {\it Let $h$ be defined as above with
$\lambda = 3/8$ and let $X(\cdot )$ be defined by (\ref{1.17e}). Let
$u_a$, $A_a$, $R_1$ and $R_2$ be sufficiently regular (for the following estimates to make sense) and satisfy the estimates
\beq
\label{2.17e}
\parallel \partial_t^j \nabla^k u_a(t) \parallel_r \ \leq \ c\ t^{-\delta (r)} \quad \hbox{\it for $2 \leq r \leq \infty$}
\eeq

\noi and in particular
\beq
\label{2.18e}
\parallel u_a\parallel_3\ \leq c_3 \ t^{-1/2}\quad , \quad \parallel \nabla u_a\parallel_4\ \leq c_4 \ t^{-3/4}\ ,
\eeq
\beq
\label{2.19e}
\parallel \nabla^2 u_a(t) \parallel_{4} \ \vee \ \parallel \partial_t \nabla u_a (t) \parallel_{4}\ \leq c\ t^{-3/4}\ ,
\eeq
\beq
\label{2.20e}
\parallel \partial_t^j \nabla^k A_a(t) \parallel_{\infty} \ \leq a\ t^{-1}\ ,
\eeq
\beq
\label{2.21e}
\parallel \partial_t^j \nabla^k R_1;L^1([t , \infty ),L^2) \parallel \ \leq r_1\ h(t)\ ,
\eeq
\beq
\label{2.22e}
\parallel R_2;L^{4/3}([t, \infty ), W_{4/3}^1) \parallel\  \leq r_2\ h(t)\ ,
\eeq

\noi for $0 \leq j+k \leq 1$, for some constants $c$, $c_3$, $c_4$, $a$, $r_1$ and $r_2$ with $c_3$ and $c_4$
sufficiently small and for all $t \geq T_0 \geq 1$. Then there exists $T$, $T_0 \leq T < \infty$
and there exists a unique solution $(v, B)$ of the
system (\ref{1.4e}) in $X([T, \infty ))$. If in addition 
\beq
\label{2.23e}
\parallel R_2;L^{1}([T, \infty ), L^{2}) \parallel\ \leq r_2\ t^{-1/2}\ h(t)\ ,
\eeq

\noi then $\nabla B$, $\partial_t B \in {\cal C}([T, \infty ), L^2)$ and $B$ satisfies the estimate 
\beq
\label{2.24e}
\parallel \nabla B(t)\parallel_2\ \vee \ \parallel \partial_t B(t)\parallel_2\ \leq C \left (  t^{-1/2}  + t^{1/4}\ h(t) \right )  h(t)
\eeq

\noi for some constant $C$ and for all $t \geq T$.}\\

\noi {\bf Proof.} We follow the sketch given in the introduction. Let
$T_0 \leq T < \infty$ and let $(v,B)\in X([T, \infty ))$. In particular
$(v, B)$ satisfies
\beq
\label{2.25e}
\parallel v(t) \parallel_2 \ \leq N_0\ h(t)
\eeq
\beq
\label{2.26e}
\parallel v ; L^{4}(J, L^3) \parallel \ \vee \ \parallel v ; L^{8/3}(J, L^4) \parallel\ \leq N_1\ h(t)
\eeq
\beq
\label{2.27e}
\parallel B ; L^{4}(J, L^4) \parallel \ \leq N_2\ h(t)
\eeq
\beq
\label{2.28e}
\parallel \partial_t v(t) \parallel_2 \ \leq N_3\ h(t)
\eeq
\beq
\label{2.29e}
\parallel \nabla v ; L^{4}(J, L^3) \parallel \ \vee \ \parallel \nabla v ; L^{8/3}(J, L^4) \parallel\ \leq N_4\ h(t)
\eeq
\beq
\label{2.30e}
\parallel \Delta v(t) \parallel_2 \ \leq N_5\ h(t)
\eeq
\beq
\label{2.31e}
\parallel \nabla B ; L^{4}(J, L^4) \parallel \ \vee \ \parallel \partial_t B  ; L^{4}(J, L^4) \parallel\ \leq N_6\ h(t)
\eeq

\noi for some constants $N_i$, $0 \leq i \leq 6$ and for all $t \geq T_0$, with $J = [t, \infty )$. Furthermore from (\ref{2.25e}) (\ref{2.30e}) it follows that 
\beq
\label{2.32e}
\parallel \nabla v(t) \parallel_2 \ \leq \left ( N_0N_5\right )^{1/2} \ h(t) \equiv N_{1/2}\ h(t) 
\eeq

\noi for all $t \geq T_0$. We first
construct a solution $(v',B')$ of the system (\ref{1.7e}) in $X([T,
\infty ))$. For that purpose, we take $t_0$, $T < t_0 < \infty$ and we
solve the system (\ref{1.7e}) in $X(I)$ where $I = [T, t_0]$ with
initial condition $(v', B')(t_0) = 0$. Let $(v'_{t_0}, B'_{t_0})$ be
the solution thereby obtained. The existence of $v'_{t_0}$ follows from
Proposition 2.1 with $V = g(|u|^2)$ and $f = G_1 - R_1$. We want to
take the limit of $(v'_{t_0}, B'_{t_0})$ as $t_0 \to \infty$ and for
that purpose we need estimates of $(v'_{t_0}, B'_{t_0})$ in $X(I)$ that
are uniform in $t_0$. Omitting the subscript $t_0$ for brevity we
define 
\beq
\label{2.33e}
N'_0 =  \ \mathrel{\mathop {\rm Sup}_{t\in I}}\ h(t)^{-1} \parallel v'(t) \parallel_2
\eeq 
\beq
\label{2.34e}
N'_1 =  \ \mathrel{\mathop {\rm Sup}_{t\in I}}\ h(t)^{-1} \left ( \parallel v';L^{4}(J,L^3) \parallel\ \vee \ \parallel v';L^{8/3}(J,L^4) \parallel\right )
\eeq 
\beq
\label{2.35e}
N'_2 =  \ \mathrel{\mathop {\rm Sup}_{t\in I}}\ h(t)^{-1} \parallel B';L^{4}(J,L^4) \parallel
\eeq
\beq
\label{2.36e}
N'_3 =  \ \mathrel{\mathop {\rm Sup}_{t\in I}}\ h(t)^{-1} \parallel \partial_t v'(t) \parallel_2
\eeq
\beq
\label{2.37e}
N'_4 =  \ \mathrel{\mathop {\rm Sup}_{t\in I}}\ h(t)^{-1} \left ( \parallel \nabla v';L^{4}(J,L^3) \parallel\ \vee \ \parallel \nabla v';L^{8/3}(J,L^4) \parallel\right )
\eeq
\beq
\label{2.38e}
N'_5 =  \ \mathrel{\mathop {\rm Sup}_{t\in I}}\ h(t)^{-1} \parallel \Delta v'(t) \parallel_2
\eeq
\beq
\label{2.39e}
N'_6 =  \ \mathrel{\mathop {\rm Sup}_{t\in I}}\ h(t)^{-1} \left ( \parallel \nabla B';L^{4}(J,L^4) \parallel\ \vee \ \parallel \partial_t B';L^{4}(J,L^4) \parallel\right )
\eeq    

\noi where $J = [t, \infty ) \cap I$ and we set out to estimate the various $N'_i$. We also define the auxiliary quantities
\bea
\label{2.40e}
&&N'_{1/2} = \ \mathrel{\mathop {\rm Sup}_{t\in I}}\ h(t)^{-1} \parallel \nabla v'(t) \parallel_2 \\
&&\widetilde{N}'_{1/2} = \ \mathrel{\mathop {\rm Sup}_{t\in I}}\ h(t)^{-1} \parallel \nabla_A v'(t) \parallel_2 
\label{2.41e}
\eea

\noi so that in particular $N'_{1/2} \leq (N'_0 N'_5)^{1/2}$. \par

We shall use the notation
$$\parallel f;L^q(J,L^r)\parallel\ = \ \parallel \ \parallel f\parallel_r;L^q(J)\parallel\ = \ \parallel \ \parallel f \parallel_r \ \parallel_q \ ,$$

\noi namely with the inner norm taken in $L^r({I\hskip-1truemm R}^3)$
and the outer norm taken in $L^q(J)$. Furthermore we shall use a shorthand
notation for two important cases, namely
$$\parallel \cdot\  ; L^1(J,L^2)\parallel\ = \ \parallel \cdot \parallel_+ \quad \hbox{and} \quad \parallel\cdot \ ; L^{4/3}(J,L^{4/3})\parallel\ = \ \parallel \cdot \parallel_*\ .$$

\noi We first estimate $N'_0$, defined by (\ref{2.33e}). From (\ref{2.2e}) we obtain 
$$\parallel v'(t)\parallel_2\ \leq \ \parallel G_1 \parallel_+ \ + \ \parallel R_1 \parallel_+$$

\noi with $G_1$ defined by (\ref{1.5e}). We estimate
$$\parallel B \cdot \nabla u_a \parallel_+ \ \leq\ \parallel \ \parallel B \parallel_4 \ \parallel\nabla u_a \parallel_4 \ \parallel_1 \ \leq C\ c_4\ N_2\ h(t)$$

\noi by Lemma 2.2,
$$\parallel B \cdot A_a u_a \parallel_+ \ \leq \ \parallel\ \parallel B \parallel_4\ \parallel A_a \parallel_{\infty} \ \parallel u_a \parallel_4 \ \parallel_1 $$
$$\leq c\ a\ N_2\ h \parallel t^{-7/4}\parallel_{4/3} \ \leq c\ a \ N_2\ t^{-1}\ h(t) \ ,$$
$$\parallel B^2 u_a \parallel_+ \ \leq \ \parallel\ \parallel B \parallel_4^2\ \parallel u_a \parallel_{\infty} \ \parallel_1\ \leq c\ N_2^2\ h^2 \parallel t^{-3/2} \parallel_2 \ \leq c\ N_2^2\ t^{-1}\ h(t)^2 \ , $$

$$\parallel g(\overline{u}_av)u_a \parallel_+ \ \leq \ C\parallel\ \parallel v \parallel_2\ \parallel u_a \parallel_3^2 \ \parallel_1 \ \leq C\ c_3^2 \ N_0\ h(t)$$

\noi by Lemma 2.3, part (1) and Lemma 2.2,
\begin{eqnarray*}
\parallel g(|v|^2)u_a \parallel_+ &\leq& C \parallel \ \parallel v \parallel_3\ \parallel v \parallel_2\ \parallel u_a \parallel_3 \ \parallel_1\\
&\leq& C\ c_3\ N_0\ N_1\ t^{1/4}\ h(t)^2
\end{eqnarray*}

\noi by Lemma 2.3, part (1) and Lemma 2.2 again.\par

Collecting the previous estimates yields
\bea
\label{2.42e}
&&N'_0 \leq C_0 \Big ( c_4\ N_2 + c\ a\ N_2\ T^{-1} + c_3^2\ N_0 + c \ N_2^2\ T^{-1}\ h(T) \nn \\
&&+ c\ N_0\ N_1\ T^{-1/8}\ \overline{h}(T) + r_1 \Big )
\eea

\noi which is of the form 
\beq
\label{2.43e}
N'_0 \leq C_0 \left ( c_4\ N_2 + c_3^2\ N_0 + r_1 + (o(1);N_0, N_1, N_2)\right )
\eeq

\noi where $(o(1);\cdot , \cdots , \cdot )$ denotes a quantity depending on the variables indicated and tending to zero as $T \to \infty$ when those variables are fixed. \par

We next estimate the Strichartz norms of $v'$, namely $N'_1$ defined by (\ref{2.34e}). By Lemma 2.1, in addition to the contribution of $G_1 - R_1$ estimated above, we need to estimate
$$iA \cdot \nabla v' + (1/2) A^2 v' + g(|u|^2)v'$$

\noi in some $L^{\overline{q}}(J, L^{\overline{r}})$ for admissible $(q, r)$. We estimate
$$\parallel A_a \cdot \nabla v'\parallel_+ \   \leq \  a \ N'_{1/2} \parallel t^{-1} h\parallel_1 \ \leq 3a\ N'_{1/2}\ h(t) \ ,$$
$$\parallel A_a^2 v'\parallel_+ \   \leq \  a^2\  N'_{0} \parallel t^{-2} h\parallel_1 \ \leq a^2\ N'_{0}\ t^{-1}\ h(t) \ ,$$
$$\parallel B \cdot \nabla v';L^{8/5}(J, L^{4/3})\parallel \   \leq \  \parallel \ \parallel B \parallel_{4}\  \parallel \nabla v' \parallel_2 \ \parallel_{8/5} \ \leq C\ N_2\ N'_{1/2}\ h(t) \ \overline{h}(t) \ ,$$

\noi by Lemma 2.2,
$$\parallel B^2  v';L^2(J, L^{6/5})\parallel \   \leq \  C\parallel\ \parallel B \parallel_4^{3/2}\  \parallel \nabla B\parallel_4^{1/2}\ \parallel v' \parallel_2\ \parallel_2 \ \leq C\ N_2^{3/2} \ N_6^{1/2} \ N'_0\ h(t)^3$$

\noi by Sobolev inequalities and Lemma 2.2,
$$\parallel g(|u_a|^2)v'\parallel_+ \   \leq \  \parallel \  \parallel g(|u_a|^2)\parallel_{\infty} \ \parallel v'\parallel_2\ \parallel_1\ \leq C\ c^2 \ N'_0 \ h(t)$$

\noi by Lemma 2.3, part (2), 
$$\parallel g(|v|^2)v';L^{4/3}(J,L^{3/2})\parallel \   \leq \  C\parallel \  \parallel v\parallel_{3} \ \parallel v\parallel_2\ \parallel v' \parallel_2 \ \parallel_{4/3}\  \leq C\ N_0 \ N_1\ N'_0 \ t^{1/2}\ h(t)^3$$

\noi by Lemma 2.3, part (1) and Lemma 2.2. The term $g(\overline{u}_av)v'$ need not be considered because it is controlled by the previous ones.\par

Collecting the previous estimates yields
$$N'_1 \leq C_1\Big \{ c_4\ N_2 + c\ a\ N_2 \ T^{-1} + c_3^2\ N_0 + c\ N_2^2 \ T^{-1}\ h(T) + \ c\ N_0\ N_1 \ T^{-1/8} \ \overline{h}(T)$$
$$ + \left ( c^2 + a^2 \ T^{-1} + N_2^{3/2} \ N_6^{1/2} \ h(T)^2 + N_0\ N_1\ T^{-1/4} \ \overline{h}(T)^2 \right ) N'_0 + r_1$$
\beq
\label{2.44e}
+ \ a\ N'_{1/2} + N_2\ N'_{1/2} \ \overline{h}(T) \Big \}
\eeq

\noi which is of the form 
\beq
\label{2.45e}
N'_1 \leq C_1 \left ( c_4 N_2 + c_3^2 N_0 + c^2 N'_0 + a N'_{1/2} + r_1 +  (o(1);N_0,N_1,N_2,N_6,N'_0,N'_{1/2}) \right ) \ .
\eeq
We now turn to the estimates of $B'$. We first estimate $B'$ in
$L^4(J,L^4)$, namely we estimate $N'_2$ defined by (\ref{2.35e}), by the use of (\ref{1.5e}) (\ref{2.8e}). For that purpose
we estimate $G_2$ in $L^{4/3}(J,L^{4/3})$. The linear terms in $v$ are
estimated by 
\begin{eqnarray*}
\parallel \overline{v} \nabla u_a\parallel_* &\leq& \parallel \ \parallel v \parallel_2\  \parallel \nabla u_a \parallel_4 \ \parallel_{4/3} \\
&\leq& c_4\ N_0 \parallel t^{-3/4} h \parallel_{4/3} \ \leq \ 2c_4\ N_0\ h(t) \ ,
\end{eqnarray*}
\begin{eqnarray*}
\parallel \overline{v} A_a u_a\parallel_* &\leq& \parallel \ \parallel v \parallel_2\  \parallel A_a \parallel_{\infty}\  \parallel u_a \parallel_4\ \parallel_{4/3} \\
&\leq& a\ c\ \ N_0 \parallel t^{-7/4} h \parallel_{4/3} \ \leq \ a\ c\ N_0\ t^{-1}\ h(t) \ .
\end{eqnarray*}

\noi The linear term in $B$ is estimated by
\begin{eqnarray*}
\parallel B| u_a|^2 \parallel_* &\leq& \parallel \ \parallel B \parallel_4\  \parallel u_a \parallel_4^2\  \parallel_{4/3} \\
&\leq& c^2\ \ N_2\ h  \parallel t^{-3/2}\parallel_{2} \ \leq \ c^2\ N_2\ t^{-1}\ h(t) \ .
\end{eqnarray*}

\noi The quadratic terms in $v^2$ are estimated by
$$\parallel \overline{v} \nabla v \parallel_* \leq \parallel \ \parallel v \parallel_4\  \parallel \nabla v \parallel_{2}\   \parallel_{4/3} \ \leq   C\ N_1\ N_{1/2}\ h(t)\ \overline{h}(t)\ ,$$

\noi by Lemma 2.2, 
\begin{eqnarray*}
\parallel A_a| v|^2 \parallel_* &\leq& \parallel \ \parallel v \parallel_4\  \parallel v \parallel_2\  \parallel A_a \parallel_{\infty}\ \parallel_{4/3} \\
&\leq& a\ N_0\ N_1\ h  \parallel t^{-1}h\parallel_{8/3} \ \leq \ a\ N_0\ N_1\ t^{-5/8}\ h(t)^2 \ .
\end{eqnarray*}

\noi The quadratic terms in $Bv$ need not be considered because
$$2|\overline{v}Bu_a| \leq \left | B|u_a|^2\right | + \left | B| v|^2\right | \ .$$

\noi The cubic term $B|v|^2$ is estimated by
$$\parallel B|v|^2\parallel_*\  \leq \  \parallel B;L^4(L^4) \parallel\ \parallel v;L^3(L^r)\parallel^{3/2} \ \parallel v;L^{\infty}(L^6)\parallel^{1/2}$$
$$\leq C\ N_2\ N_1^{3/2}\ N_{1/2}^{1/2}\ h(t)^3$$

\noi where $3 < r = 18/5 < 4$ so that $(3,r)$ is an admissible pair and that the middle norm is controlled by $N_1$. \par

Collecting the previous estimates yields
\bea \label{2.46e}
&& N'_2 \leq C_2 \Big \{ c_4\ N_0 + a \ c\ N_0 \ T^{-1} + c^2\ N_2 \ T^{-1} + r_2\nn \\
&&+ N_1\ N_{1/2}\ \overline{h}(T) + a\  N_0\ N_1 \ T^{-5/8} \ h(T) + N_2\ N_1^{3/2} \ N_{1/2}^{1/2}\ h(T)^2 \Big \}
\eea

\noi which is of the form
\beq \label{2.47e}
N'_2 \leq C_2 \left ( c_4\ N_0 + r_2 +  ( o(1);N_0,N_1,N_{1/2}, N_2) \right ) \ .
\eeq

\noi We next complete the estimates of $B'$ by estimating $\nabla
B'$ and $\partial_t B'$ in $L^4(J, L^4)$, namely we estimate $N'_6$ defined by (\ref{2.39e}), through the use of (\ref{1.5e})
(\ref{2.9e}). For that purpose we estimate $\nabla G_2$ in
$L^{4/3}(J,L^{4/3})$. Now

\bea \label{2.48e}
&&\nabla G_2 = 2 \ P \ {\rm Im} \Big ( (\nabla \overline{v})\nabla_Av + (\nabla \overline{v}) \nabla_A u_a + (\nabla \overline{u}_a) \nabla_A v \Big ) \nn \\
&&- P (\nabla A) \left ( |v|^2 + 2\ {\rm Re}\ \overline{u}_a v\right ) - P (\nabla B) |u_a|^2 - 2\ P\ B \ {\rm Re}\ \overline{u}_a \nabla u_a \ .
\eea

\noi The estimate of $\nabla G_2$ in $L^{4/3}(J, L^{4/3})$ proceeds
exactly as that of $G_2$ in the same space, with one additional gradient
acting on each factor in each term, except for two facts. First because
of the symmetry of the quadratic form $P$ Im $(\overline{v}_1 \nabla_A
v_2)$, we can always ensure that no terms occur with two derivatives on
$v$ or $u_a$. Second, the quadratic terms coming from $\overline{v}
Bu_a$ have to be estimated explicitly because they are no longer
estimated by polarisation. When hitting $v$, and additional gradient
produces a replacement of $N_0$ by $N_{1/2}$ and of $N_1$ by $N_4$ in
the estimates. When hitting $B$, it produces a replacement of $N_2$ by
$N_6$. When hitting $u_a$ or $A_a$, it only requires higher regularity
of these functions, but does not change the form of the estimates. With
those remarks available, only the terms from $\nabla
(\overline{v}Bu_a)$ and from $B\nabla |v|^2$ need new estimates.\par

The linear terms in $v$ are estimated by 
\begin{eqnarray*}
&&\parallel ( \nabla \overline{v})\nabla u_a\parallel_* \ \leq 2 c_4 \ N_{1/2}\ h(t) \ ,\\
&&\parallel ( \nabla \overline{v}) A_a u_a\parallel_* \ \leq a\ c \ N_{1/2}\ t^{-1}\ h(t) \ ,\\
&&\parallel \overline{v}( \nabla A_a) u_a + \overline{v} A_a\nabla u_a \parallel_* \ \leq 2a\ c \ N_{0}\ t^{-1}\ h(t) \ .
\end{eqnarray*} 

The linear terms in $B$ are estimated by
\begin{eqnarray*}
&&\parallel ( \nabla B)|u_a|^2\parallel_* \ \leq c^2 \ N_{6}\ t^{-1}\ h(t) \ ,\\
&&\parallel B \overline{u}_a \nabla u_a\parallel_* \ \leq c^2 \ N_{2}\ t^{-1}\ h(t) \ .
\end{eqnarray*}

The quadratic terms in $v^2$ are estimated by 
\begin{eqnarray*}
&&\parallel | \nabla v|^2\parallel_* \ \leq \ \parallel \ \parallel \nabla v \parallel_4\ \parallel  \nabla v \parallel_2 \ \parallel_{4/3}\ \leq C \ N_4\ N_{1/2}\ h(t) \ \overline{h}(t)\ , \\
&&\parallel ( \nabla \overline{v})A_a v\parallel_* \ \leq \ \parallel \ \parallel \nabla v \parallel_2\ \parallel  A_a \parallel_{\infty}  \  \parallel v \parallel_4\ \parallel_{4/3}\ \leq \ a\ N_{1/2}\ N_1 \ t^{-5/8} \ h(t)^2  \ ,\\
&&\parallel ( \nabla A_a)|v|^2\parallel_* \ \leq \ a\ N_0\ N_1 \ t^{-5/8} \ h(t)^2  \ .
\end{eqnarray*}

The quadratic terms in $Bv$ are estimated by 
\begin{eqnarray*}
&&\parallel ( \nabla \overline{v})Bu_a\parallel_* \ \leq \ \parallel \ \parallel \nabla v \parallel_2\ \parallel  B \parallel_4 \ \parallel  u_a \parallel_{\infty}\ \parallel_{4/3}\\
&&\leq \ c\ N_{1/2}\ N_2 \ h \parallel  t^{-3/2} h \parallel_2 \ \leq c\ N_{1/2}\ N_2 \ t^{-1} \ h(t)^2  
\end{eqnarray*}  

\noi and similarly
\begin{eqnarray*}
&&\parallel  \overline{v} (\nabla B)u_a\parallel_* \ \leq c \ N_{0}\ N_6 \ t^{-1}\ h(t)^2 \ , \\
&&\parallel  \overline{v} B \nabla u_a\parallel_* \ \leq c \ N_{0}\ N_2 \ t^{-1}\ h(t)^2 \ .
\end{eqnarray*}

The cubic terms from $B|v|^2$ are estimated by
$$\parallel ( \nabla B)|v|^2\parallel_* \ \leq C \ N_{6}\ N_1^{3/2}\ N_{1/2}^{1/2}\ h(t)^3 \ ,$$
\begin{eqnarray*}
\parallel B \overline{v} \nabla v\parallel_* &\leq & C \parallel B;L^4(L^4)\parallel \ \parallel v;L^4(L^3)\parallel^{1/2} \ \parallel \nabla v; L^4(L^3)\parallel^{3/2}\\
  &\leq &C \ N_2\ N_{1}^{1/2}\ N_4^{3/2}\ h(t)^3 \ .
\end{eqnarray*}

Collecting the previous estimates yields
\bea \label{2.49e}
N'_6 &\leq& C_6 \Big \{ c_4\ N_{1/2} + a c (N_{1/2} + N_0) T^{-1} + c^2 (N_2 + N_6) T^{-1} + r_2 + N_{1/2} \ N_4 \ \overline{h}(T) \nn \\
&&+ a (N_{1/2} + N_0) N_1\ T^{-5/8} h(T) + c \left ( N_2\ N_{1/2} + N_2\ N_0 + N_6 \ N_0 \right )  T^{-1}\ h(T)\nn \\
&& + \left ( N_6\ N_1^{3/2} \ N_{1/2}^{1/2} + N_2 \ N_1^{1/2} \ N_4^{3/2}\right ) h(T)^2 \Big \}\eea

\noi which is of the form
\beq \label{2.50e}
N'_6 \leq C_6 \left ( c_4\ N_{1/2} + r_2 + (o(1); N_0, N_1, N_{1/2}, N_4, N_2, N_6) \right ) \ .
\eeq

We now come back to the estimates of $v'$ and we first estimate
$\partial_tv'$ in $L^2$, namely we estimate $N'_3$ defined by
(\ref{2.36e}) by using (\ref{2.3e}). Here however we encounter a
technical difficulty due to the fact that $B$ a priori does not satisfy
the assumption $\nabla B \in L_{loc}^1(I, L^{\infty})$ needed in
Proposition 2.1, part (2) in order to derive (\ref{2.3e}). We
circumvent that difficulty by first regularizing $B$, introducing the
associated solution $v'$ which then satisfies (\ref{2.3e}), deriving the
$N'_3$ estimate for the auxiliary solution, and removing the
regularization by a limiting procedure, which preserves the estimate.
Here in order not to burden the proof with technicalities, we provide
only the derivation of the estimates from (\ref{2.3e}) and we refer to
the proof of Proposition 3.2, part (1) in \cite{7r} for the technical
details. From (\ref{2.3e}) (\ref{2.4e}) with $V = g(|u|^2)$ and $f =
G_1 - R_1$, we obtain 
\beq \label{2.51e}
\parallel  \partial_t v' \parallel _2\ \leq \ \parallel i(\partial_t A)\cdot \nabla_Av' + \left ( \partial_t g(|u|^2)\right ) v' + \partial_t G_1 - \partial_t R_1 \parallel_+\ + \ \parallel (G_1 - R_1)(t_0)\parallel_2
\eeq

\noi with $G_1$ defined in (\ref{1.5e}).\par

We first estimate the terms containing $v'$, starting with $i(\partial_tA)\cdot \nabla_A v'$.
\begin{eqnarray*}
\parallel (\partial_t A_a) \cdot \nabla_A v'\parallel_+ &\leq & \parallel\  \parallel \partial_t A_a\parallel_{\infty} \ \parallel \nabla_A v'\parallel_2 \ \parallel_1\\
  &\leq &a \ \widetilde{N}'_{1/2}\parallel t^{-1} h \parallel_1 \ \leq \ 3 a\ \widetilde{N}'_{1/2}\ h(t) \ ,
\end{eqnarray*}
$$\parallel (\partial_t B) \cdot \nabla v'\parallel_+ \  \leq \  \parallel\  \parallel \partial_t B\parallel_{4} \ \parallel \nabla v'\parallel_4 \ \parallel_1\ \leq C\ N_6\ N'_4 \ h(t) \ \overline{h}(t) \ ,$$
\begin{eqnarray*}
\parallel (\partial_t B) \cdot A_a v'\parallel_+ & \leq &  \parallel\  \parallel \partial_t B\parallel_{4} \ \parallel A_a\parallel_{\infty} \ \parallel v'\parallel_4\ \parallel_1\\
  &\leq& a\ N_6 \ N'_1 \ h^2  \parallel t^{-1} \parallel_{8/3} \ \leq a\ N_6\ N'_1 \ t^{-5/8}\ h(t)^2  \ ,
\end{eqnarray*}
\begin{eqnarray*}
\parallel (\partial_t B) \cdot B v'\parallel_+ & \leq  & C\parallel\  \parallel \partial_t B\parallel_{4} \left ( \parallel B\parallel_4 \ \parallel \nabla B \parallel_4^3 \right)^{1/4}\  \parallel v' \parallel_4 \ \parallel_1\\
  &\leq& C\ N_6^{7/4} \ N_2^{1/4} \ N'_1\ t^{1/8} \ h(t)^3 
\end{eqnarray*}

\noi by Lemma 2.2.\par

We next estimate the terms coming from $(\partial_tg(|u|^2))v'$.
$$\parallel g\left ( \overline{u}_a \partial_t u_a\right ) v' \parallel_+ \ \leq \  \parallel\ \parallel g\left ( \overline{u}_a \partial_t u_a\right )\parallel_{\infty} \ \parallel v'\parallel_2\ \parallel_1 \ \leq C\ c^2\ N'_0\ h(t)\ ,$$
$$\parallel g\left ( \overline{u}_a \partial_t v\right ) v' \parallel_+ \ \leq \ \parallel\ \parallel g\left ( \overline{u}_a \partial_t v\right )\parallel_{\infty} \ \parallel v'\parallel_2\ \parallel_1\ \leq C\ c\ N_3\ N'_0\ h(t)^2\ ,$$
$$\parallel g\left ( (\partial_t \overline{u}_a)  v\right ) v' \parallel_+ \ \leq \ C\ c\ N_0\ N'_0\ h(t)^2$$

\noi by Lemma 2.3, part (2) and Lemma 2.2,
$$\parallel g\left ( \overline{v} \partial_t v\right ) v' \parallel_+ \ \leq \ \parallel\ \parallel v\parallel_{3} \ \parallel \partial_t v\parallel_2\ \parallel v'\parallel_3\ \parallel_1\ \leq\  C\ N_1\ N_3\ N'_1\ t^{1/2}\ h(t)^3$$

\noi by Lemma 2.3, part (1) and Lemma 2.2.\par

We next estimate $\partial_tG_1$. The estimates are similar to those
performed when estimating $v'$ in $L^2$, with an additional time
derivative acting on each factor in each term. This has the effect of
requiring more regularity on $(A_a, u_a)$ when that derivative hits
$(A_a, u_a)$, without changing the form the estimate, and of replacing
one factor $N_2$ by $N_6$ when that derivative hits $B$ and one factor $N_0$ by $N_3$ when that derivative hits $v$. Thus we obtain
$$\parallel (\partial_t B) \cdot \nabla_{A_a} u_a \parallel_+ \ \leq N_6 (C\ c_4 + a\ c \ t^{-1}) h(t) \ ,$$
$$\parallel B \cdot \partial_t  \nabla_{A_a} u_a \parallel_+ \ \leq c\ N_2 (C+ a \ t^{-1}) h(t) \ ,$$
$$\parallel (\partial_t B) B u_a \parallel_+ \ \leq c\ N_2\ N_6  \ t^{-1}\  h(t)^2 \ ,$$
$$\parallel B^2 \partial_t   u_a \parallel_+ \ \leq c\ N_2^2  \ t^{-1} \ h(t)^2 \ .$$
$$\parallel \partial_t \left ( g(\overline{u}_av) u_a\right ) \parallel_+ \ \leq \ C\ c_3\left ( c_3 N_3 + c\ N_0 \right ) h(t)$$
$$\parallel \partial_t \left ( g(|v|^2) u_a \right )\parallel_+ \ \leq \ C\ N_1 \left ( c_3\ N_3 + c\ N_0\right ) t^{1/4}\ h(t)^2\ .$$

\noi We finally estimate $\parallel \partial_tv'(t_0)\parallel_2$ and for that
purpose we need pointwise (in time) estimates of $R_1$ and of $B$. Now
from (\ref{2.21e}) it follows that
\beq
\label{2.52e}
\parallel R_1(t) \parallel_2\ \leq \ \parallel \partial_t R_1 \parallel_+ \ \leq \ r_1 \ h(t)
\eeq

\noi while from (\ref{2.27e}) (\ref{2.31e})
$$\parallel B(t)\parallel_4^4 \ \leq \  4 \int_t^{\infty} dt'\parallel B(t')  \parallel_4^3 \ \parallel \partial_t B(t') \parallel_4 \ \leq 4N_2^3\ N_6 \ h(t)^4$$

\noi and therefore
$$\parallel B(t)\parallel_4 \ \leq  \widetilde{N}_2\ h(t) \equiv \sqrt{2} \left ( N_2^3\ N_6\right )^{1/4}\ h(t) \ .$$

\noi We then estimate

$$\parallel G_1\parallel_2 \ \leq \  \parallel B\parallel_{4} \left (  \parallel \nabla u_a\parallel_4 \ + \ \parallel A_a\parallel_{\infty} \ \parallel u_a\parallel_4\right ) + \ \parallel B \parallel_4^2\ \parallel u_a\parallel_{\infty}$$
$$+ \ \parallel g\left ( |v|^2 + 2{\rm Re}\ \overline{u}_av\right ) \parallel_6 \ \parallel u_a\parallel_3$$
$$\leq c\ \widetilde{N}_2 (1 + a \ t^{-1})t^{-3/4} h(t) + c\ \widetilde{N}_2^2 \ t^{-3/2} \ h(t)^2$$
\beq
\label{2.53e}
+ C\Big ( c_3^2\ N_0 \ t^{-1}\ h(t) + c_3\ N_0^{3/2}\ N_{1/2}^{1/2}\ t^{-1/2}\ h(t)^2\Big ) 
\eeq

\noi by Lemma 2.3 for the terms containing $g$ and the definitions.\par

Collecting the previous estimates and in particular (\ref{2.52e}) (\ref{2.53e}) taken at $t_0 \geq t$, we obtain
$$N'_3 \leq C_3 \Big \{ a\ \widetilde{N}'_{1/2} + c_4\ N_6 + c\ N_2 + c\ a (N_6 + N_2)T^{-1} + c_3^2\ N_3 + c^2(N_0 + N'_0) + r_1$$
$$+ N_6\ N'_4 \ \overline{h}(T) + a\ N_6\ N'_1 \ T^{-5/8} h(T) + c\ N_2 (N_6 + N_2) T^{-1} h(T) $$
$$+ c\left ( N_3 + N_0\right ) N'_0 \ h(T) + c\left ( N_3 + N_0\right ) N_1\ T^{-1/8}\ \overline{h}(T)$$
$$+ N_6^{7/4}\ N_2^{1/4} N'_1\ T^{-1/4} \ h(T) \ \overline{h}(T) + N_1\ N_3 \ N'_1 \ T^{-1/4} \ \overline{h}(T)^2$$
$$+ c \ \widetilde{N}_2 (1 + a \ T^{-1})T^{-3/4} + c\ \widetilde{N}_2^2\ T^{-3/2}\ h(T)$$
\beq \label{2.54e}
+ c_3^2\ N_0 \ T^{-1} + c_3 \ N_0^{3/2}\ N_{1/2}^{1/2} \ T^{-1/2}\ h(t) \Big \}
\eeq

\noi which is of the form 
\bea \label{2.55e}
N'_3 &\leq &C_3 \Big ( a \widetilde{N}'_{1/2} + c_4 \ N_6 + c\ N_2 + c_3^2\ N_3 + c^2\left ( N_0 + N'_0\right ) + r_1 \nn\\
&&+ ( o(1);N_0,N'_0,N_1,N'_1,N_{1/2}, N_2, \widetilde{N_2},N_3,N'_4,N_6 )\Big ) \ .
\eea

We next estimate $\parallel \Delta_Av'\parallel_2$\ . From (\ref{1.7e}), we obtain
\begin{eqnarray*}
\parallel \Delta_A v' \parallel &\leq& 2 \left ( \parallel \partial_t v' \parallel_2\ + \ \parallel g(|u|^2)v' \parallel_2\ + \ \parallel G_1 \parallel_2 \ + \ \parallel R_1 \parallel_2 \right )\\
&\leq& 2\left ( N'_3 + r_1 \right ) h(t) + \ 2\parallel g(|u|^2)v' \parallel_2\ + \ 2 \parallel G_1 \parallel_2\ .
\end{eqnarray*}

\noi Furthermore
$$\parallel g(|u|^2)v'\parallel_2 \ +\  \parallel G_1\parallel_2 \ \leq \  \parallel g(|u|^2)\parallel_{\infty} \ \parallel v'\parallel_2 \ + \ \parallel G_1 \parallel_2 \ $$
$$\leq c\ \widetilde{N}_2 (1 + a \ t^{-1})t^{-3/4} h(t) + c\ \widetilde{N}_2^2 \ t^{-3/2} \ h(t)^2 + C \Big ( c^2\ N'_0 \ t^{-1}\ h(t) $$
\beq
\label{2.56e}
+ N_0\ N_{1/2}\ N'_0\ h(t)^3 + c_3^2\ N_0 \ t^{-1}\ h(t) + c_3\ N_0^{3/2}\ N_{1/2}^{1/2}\ t^{-1/2}\ h(t)^2\Big ) \equiv M_1\ h(t)
\eeq

\noi by Lemma 2.3 for the terms containing $g$ and by (\ref{2.53e}), so that
\bea
\label{2.57e}
\parallel \Delta_A v'(t)\parallel_2 &\leq& 2\left ( N'_3 + r_1 + M_1 \right ) h(t)\nn \\
&=& 2 \left ( N'_3 + r_1 + \left ( o(1);\widetilde{N}_2,N_0,N'_0,N_{1/2}\right ) \right ) h(t) \ .
\eea 

\noi As a consequence,
\bea
\label{2.58e}
\widetilde{N}'_{1/2} &\leq& \left ( 2N'_0 (N'_3 + r_1 + M_1(T))\right )^{1/2}\nn \\
&\leq& \left ( 2N'_0\left ( N'_3 + r_1 + \left ( o(1);\widetilde{N}_2,N_0,N'_0,N_{1/2} \right ) \right )\right )^{1/2}\ .
\eea

\noi We next estimate $\parallel \Delta v'(t)\parallel_2$\ , namely $N'_5$ defined by (\ref{2.38e}). From 
$$\Delta_Av' = \Delta v' - 2i A_a \cdot \nabla_Av' - 2i B\cdot \nabla v' + (A_a^2 - B^2) v'$$

\noi we obtain
\begin{eqnarray*}
\parallel \Delta v'\parallel_2  &\leq&   \parallel \Delta_A v'\parallel_2 \ + \ 2 \parallel A_a\parallel_{\infty} \ \parallel \nabla_A v'\parallel_2 \ + \ \parallel A_a\parallel_{\infty}^2 \ \parallel v'\parallel_2\\
&&+ \ 2\parallel B\parallel_4  \ \parallel \nabla v'\parallel_4\ +\ \parallel B\parallel_4^2 \ \parallel v'\parallel_{\infty}\ .
\end{eqnarray*}

\noi Now
$$\parallel \nabla v'\parallel_4 \ \leq \ C \parallel \Delta v'\parallel_2^{7/8}\ \parallel v'\parallel_2^{1/8} \ ,$$
$$\parallel v' \parallel_{\infty}\ \leq C \parallel  \Delta v'\parallel_2^{3/4}\ \parallel v' \parallel_2^{1/4}\ ,$$

\noi and therefore 
\begin{eqnarray*}
\parallel \Delta v' \parallel_2 &\leq&  (1 + \varepsilon ) \parallel \Delta_A v' \parallel_2\ + \ (1 + \varepsilon^{-1} )\parallel A_a \parallel_{\infty}^2 \ \parallel v' \parallel_2 \\
&&+ \ \varepsilon \parallel \Delta v' \parallel_2 \ + \ C_{\varepsilon} \parallel B \parallel_4^8\ \parallel v' \parallel_2\ .  
\end{eqnarray*}

\noi Taking $\varepsilon = 1/3$ yields
$$\parallel \Delta v'\parallel_2 \ \leq \ 2 \parallel \Delta_Av' \parallel_2\ + \ C\left ( \parallel A_a\parallel_{\infty}^2 \ + \ \parallel B \parallel_4^8 \right ) \parallel v'\parallel_2$$

\noi and therefore by (\ref{2.38e}) (\ref{2.57e})
\beq \label{2.59e}
N'_{5} \leq 4 \left ( N'_3 + r_1 + M_1 \right ) + C\left (a^2 \ T^{-2} + \widetilde{N}_2^8\ h(T)^8 \right ) N'_0
\eeq

\noi which is of the form
\beq \label{2.60e}
N'_{5} \leq 4 \left ( N'_3 + r_1 + (o(1);\widetilde{N}_2,N_0,N'_0,N_{1/2}  )\right ) \ .
\eeq

We finally estimate the Strichartz norms of $\nabla v'$. For that
purpose, by Lemma 2.1, we have to estimate the following quantity in the sum of
spaces of the type $L^{\overline{q}}(J,L^{\overline{r}})$ for
admissible pairs $(q, r)$~: 
\bea
\label{2.61e}
Q &=& \nabla \left ( i A \cdot \nabla v' + (1/2) A^2 v' + g(|u|^2)v' + G_1 - R_1 \right )\nn\\
&=& i A \cdot \nabla^2 v' + i \nabla A \cdot \nabla_A v' + (1/2) A^2 \nabla v' + \nabla (g(|u|^2)v') + \nabla G_1 - \nabla R_1 \ .\nn \\
\eea

\noi The estimates are similar to those performed when estimating
$\parallel v'\parallel_2$ and the Strichartz norms of $v'$ (see the proof of
(\ref{2.42e}) (\ref{2.44e})), with an additional  gradient acting on each
factor in each term, thereby producing the replacement of $N_0$ by $N_{1/2}$, of $N'_0$ by
$N'_{1/2}$, of $N'_{1/2}$ by $N'_5$
and of $N_2$ by $N_6$ at suitable places. More precisely, the terms
containing $v'$ are estimated by
$$\parallel A_a \cdot \nabla^2 v' \parallel_+ \ \leq \ \parallel \ \parallel A_a \parallel_{\infty} \ \parallel \Delta v' \parallel_2\ \parallel_1\ \leq 3 a\ N'_5\  h(t) \ ,$$
$$\parallel \nabla A_a \cdot \nabla_A v' \parallel_+ \ \leq \ \parallel \ \parallel \nabla A_a \parallel_{\infty} \ \parallel \nabla_A  v' \parallel_2\ \parallel_1\ \leq 3a\ \widetilde{N}'_{1/2}\  h(t) \ ,$$
$$\parallel A_a^2 \nabla v' \parallel_+ \ \leq \ \parallel \ \parallel A_a \parallel_{\infty}^2 \ \parallel \nabla v' \parallel_2\ \parallel_1\ \leq a^2\ N'_{1/2} \ t^{-1}\ h(t) \ ,$$
$$\parallel B \cdot \nabla^2 v';L^{8/5}(J,L^{4/3})\parallel \ \leq\  \parallel \ \parallel B \parallel_4 \ \parallel \Delta v' \parallel_2\ \parallel_{8/5}\ \leq\  C\ N_2\ N'_5 \ h(t) \ \overline{h}(t) ,$$
\begin{eqnarray*}
\parallel \nabla B \cdot \nabla_A v';L^{8/5}(J,L^{4/3})\parallel &\leq& \parallel \ \parallel \nabla B \parallel_4 \ \parallel \nabla_A v' \parallel_2\ \parallel_{8/5}\nn \\
&\leq& C\ N_6\ \widetilde{N}'_{1/2} \ h(t) \ \overline{h}(t) \ ,
\end{eqnarray*}
\begin{eqnarray*}
\parallel B^2 \nabla v';L^{2}(J,L^{6/5})\parallel &\leq& C\parallel \ \parallel B \parallel_4^{3/2} \parallel \nabla B \parallel_4^{1/2} \ \parallel \nabla v' \parallel_2\   \parallel_{2}\nn \\
&\leq& C\ N_2^{3/2}\  N_6^{1/2} \ N'_{1/2} \  h(t)^3 \ , 
\end{eqnarray*}
$$\parallel \nabla \left ( g(|u_a|^2)v'\right )\parallel_+\ \leq \ C\ c^2 \left ( N'_0 + N'_{1/2}\right ) h(t)\ ,$$
$$\parallel \nabla \left ( g(\overline{u}_av)v'\right );L^{4/3}(J,L^{3/2}) \parallel\ \leq \ C\ c \left ( N_0\ N'_0 + N_0\ N'_{1/2} + N_{1/2} \ N'_0\right ) t^{1/4}\ h(t)^2\ ,$$
$$\parallel \nabla \left ( g(|v|^2)v'\right ) ;L^{4/3}(J,L^{3/2})\parallel\ \leq \ C\ N_1 \left ( N_0\ N'_{1/2} + N_{1/2} \ N'_0\right ) t^{1/2}\ h(t)^3$$

\noi where we have used again Lemmas 2.2 and 2.3 in the estimates of the terms containing $g$.\par

The terms from $\nabla G_1$ are estimated by
$$\parallel \nabla B \cdot \nabla_{Aa} u_a\parallel_+ \ \leq \ N_6 \left ( C\ c_4 + c\ a \ t^{-1}\right ) h(t) \ ,$$
$$\parallel B \cdot \nabla \nabla_{Aa} u_a\parallel_+ \ \leq \ C\ c\ N_2 \left ( 1 + a\  t^{-1}\right ) h(t) \ ,$$
$$\parallel B^2 \nabla u_a\parallel_+ \ \leq \ c\ N_2^2 \  t^{-1}\  h(t)^2 \ ,$$
$$\parallel B ( \nabla B)u_a\parallel_+ \ \leq \ c\ N_2\ N_6  \  t^{-1}\  h(t)^2 \ .$$
$$\parallel \nabla \left ( g(\overline{u}_av)u_a\right )\parallel_+ \ \leq C\ c_3 \left ( c_3\ N_{1/2} + c\ N_{0} \right )  h(t)\ ,$$
$$\parallel \nabla \left ( g(|v|^2)u_a\right )\parallel_+\ \leq \ C\ N_1 \left ( c_3\ N_{1/2} + c\ N_0\right ) t^{1/4}\ h(t)^2\ .$$

\noi Collecting the previous estimates yields
$$N'_4 \leq C_4 \Big \{ a\left ( N'_5 + \widetilde{N}'_{1/2} \right ) + a^2 \ N'_{1/2} \ T^{-1} + c\left (  N_6 + N_2 \right ) \left ( 1 + a\ T^{-1} \right ) $$
$$+ \ c^2 \left ( N_0 + N'_0 + N_{1/2} + N'_{1/2} \right ) + r_1$$
$$+ \left ( N_2\ N'_5 + N_6 \ \widetilde{N}'_{1/2}\right ) \overline{h}(T) + c\ N_2 (N_6 + N_2)T^{-1} \ h(T)$$
$$+ \ c \left ( \left ( N_{1/2} + N_0\right ) \left ( N_1 + N'_0\right ) + N_0 \ N'_{1/2} \right ) T^{-1/8} \ \overline{h}(T)$$
\beq
\label{2.62e}
+ \ N_2^{3/2}\ N_6^{1/2}\ N'_{1/2}\ h(T)^2 + N_1 \left ( N_0\ N'_{1/2} + N_{1/2}\ N'_0 \right ) T^{-1/4} \ \overline{h}(T)^2 \Big \}
\eeq

\noi which is of the form 
\bea
\label{2.63e}
N'_4 &\leq & C_4  \Big ( a\left ( N'_5 + \widetilde{N}'_{1/2} \right ) + c (N_6 + N_2) + c^2\left ( N_0 + N'_0 + N_{1/2} + N'_{1/2}\right ) + r_1 \nn \\
&&+  ( o(1);N_0, N'_0, N_1, N_{1/2}, N'_{1/2}, \widetilde{N}'_{1/2} , N'_5, N_2 , N_6) \Big )  \ .\eea

\noi From the previous estimates, more precisely from (\ref{2.42e})
(\ref{2.44e}) (\ref{2.46e}) (\ref{2.49e}) (\ref{2.54e}) (\ref{2.59e})
(\ref{2.62e}), it follows that the $N'_i$, $0 \leq i \leq 6$ are
estimated in terms of the $N_i$, $0 \leq i \leq 6$, provided $T$ is
sufficiently large. In fact (\ref{2.42e}) (\ref{2.46e}) (\ref{2.49e})
provide estimates of $N'_0$, $N'_2$ and $N'_6$. Denoting by $C$ a
general constant depending on $T$ and on the $N_i$, it follows from
(\ref{2.44e}) (\ref{2.59e}) and from the definition of $N'_{1/2}$,
$\widetilde{N}'_{1/2}$ that 
\beq
\label{2.64e}
N'_1 \leq C\left ( 1 + N'_{1/2}\right ) \ , \ N'_{1/2} \vee \widetilde{N}'_{1/2} \leq C \left ( 1 + N'_5 \right )^{1/2}\ , \ N'_5 \leq 4 N'_3 + C
\eeq

\noi so that it remains only to estimate $N'_3$ and $N'_4$. Substituting the previous estimates into (\ref{2.54e}) (\ref{2.62e}) yields
\beq
\label{2.65e}
\left \{ \begin{array}{l} N'_3 \leq C_3 \ N_6 \ N'_4\ \overline{h}(T) + \hbox{terms sublinear in $N'_3$}\\ \\  N'_4 \leq  4 C_4 \ N'_3 \left ( a +  N_2 \ \overline{h}(T) \right ) + \hbox{terms sublinear in $N'_3$}\end{array} \right . 
\eeq

\noi which ensure the required estimate of $N'_3$, $N'_4$ provided $T$ is sufficiently large so that
\beq
\label{2.66e}
4C_3\ C_4 \left ( a + N_2\ \overline{h}(T)\right ) N_6\ \overline{h}(T) < 1\ ,\eeq

\noi which we assume from now on. Note that the terms responsible for
that large $T$ condition are the terms $\partial_tB\cdot \nabla v'$
from (\ref{2.51e}) and $A \cdot \nabla^2v'$ from (\ref{2.61e}). No such
condition was required at this stage in the simpler case of the
(WS)$_3$ system \cite{8r}. The estimates obtained for the $N'_i$ are
obviously uniform in $t_0$. \par

We now take the limit $t_0 \to \infty$ of $(v'_{t_0}, B'_{t_0})$,
restoring the subscript $t_0$ for that part of the argument. Let $T <
t_0 < t_1 < \infty$ and let $(v'_{t_0}, B'_{t_0})$ and $(v'_{t_1},
B'_{t_1})$ be the corresponding solutions of (\ref{1.7e}). From the
$L^2$ norm conservation of the difference $v'_{t_0}- v'_{t_1}$ and
from (\ref{2.42e}), it follows that for all $t \in [T, t_0]$ \beq
\label{2.67e}
\parallel v'_{t_0}(t)- v'_{t_1}(t)\parallel_2\ = \ \parallel v'_{t_1}(t_0)\parallel_2\ \leq \ K_0\ h(t_0)
\eeq 

\noi where $K_0$ is the RHS of (\ref{2.42e}), while from (\ref{1.7e}) (\ref{2.8e}) (\ref{2.9e}) (\ref{2.46e}) (\ref{2.49e}) and the initial conditions, it follows that
\bea
\label{2.68e}
&&\parallel B'_{t_0} - B'_{t_1};L^4([T, t_0], W_4^1)\parallel\ \vee \ \parallel \partial_t (B'_{t_0} - B'_{t_1});L^4([T, t_0],L^4)\parallel \nn \\
&&\leq \ C \parallel G_2 - R_2;L^{4/3}([t_0,t_1], W_{4/3}^1)\parallel\ \leq (K_2 + K_6) h(t_0)
\eea

\noi where $K_2$ and $K_6$ are the RHS of (\ref{2.46e}) and (\ref{2.49e}) respectively.\par

It follows from (\ref{2.67e}) (\ref{2.68e}) that there exists $(v', B')
\in L_{loc}^{\infty}([T, \infty ), L^2) \oplus$\break \noindent $L_{loc}^4([T, \infty ),
W_4^1)$ with $\partial_tB' \in L_{loc}^4([T, \infty ), L^4)$ such that
$(v'_{t_0}, B'_{t_0})$ converges to $(v', B')$ in that space when $t_0
\to \infty$. From the uniformity in $t_0$ of the estimates
(\ref{2.42e}) (\ref{2.46e}) (\ref{2.49e}), it follows that $(v', B')$
satisfies the same estimates in $[T, \infty )$, namely that
(\ref{2.43e}) (\ref{2.47e}) (\ref{2.50e}) hold with $N'_i$ defined by
(\ref{2.33e}) (\ref{2.35e}) (\ref{2.39e}) with $I = [T, \infty )$.
Furthermore it follows by a standard compactness  argument that $(v',
B') \in X([T, \infty ))$ and that $v'$ satisfies the remaining
estimates, namely (\ref{2.45e}) (\ref{2.55e}) (\ref{2.60e})
(\ref{2.63e}) with the remaining $N'_i$ again defined by (\ref{2.34e})
(\ref{2.36e}) (\ref{2.37e}) (\ref{2.38e}) with $I = [T, \infty )$.
Clearly $(v', B')$ satisfies the system (\ref{1.7e}).

From now on, $I$ denotes the interval $[T, \infty )$. The previous
construction defines a map $\phi : (v, B) \to (v', B')$ from $X(I)$ to
itself. The next step consists in proving that the map $\phi$ is a
contraction on a suitable closed bounded set ${\cal R}$ of $X(I)$. We
define ${\cal R}$ by the conditions (\ref{2.25e})-(\ref{2.31e}) for
some constants $N_i$ and for all $t \in I$. We first show that for a
suitable choice of $N_i$ and for sufficiently large $T$, the map $\phi$
maps ${\cal R}$ into ${\cal R}$. By (\ref{2.43e}) (\ref{2.45e})
(\ref{2.47e}) (\ref{2.50e}) (\ref{2.55e}) (\ref{2.60e}) (\ref{2.63e}), it suffices for that purpose that

\beq \label{2.69e} \left \{ \begin{array}{l} (N'_0 \leq ) \ C_0 \Big ( c_4 \ N_2 + c_3^2\ N_0 + r_1 + o(1)\Big ) \leq N_0 \\ \\ (N'_1 \leq )
  \ C_1 \Big ( c_4\ N_2 + c_3^2 \ N_0 + c^2\ N'_0 + a N'_{1/2} + r_1 +
o(1)\Big )  \leq N_1 \\ \\ (N'_2 \leq ) \ C_2 \Big ( c_4\ N_0 + r_2
+ o(1)\Big ) \leq N_2 \\ \\ (N'_6 \leq ) \ C_6 \Big ( c_4\ N_{1/2} +
r_2 + o(1) \Big ) \leq N_6 \\ \\ (N'_3 \leq ) \ C_3 \Big ( a\
\widetilde{N}'_{1/2} + c_4\ N_6 + c\ N_2 + c_3^2\ N_3 + c^2 (N_0 +
N'_1) + r_1 + o(1)\Big ) \leq N_3\\ \\ (N'_4 \leq ) \ C_4 \Big (
a\left ( N'_5 + \widetilde{N}'_{1/2}\right ) + c(N_6 + N_2) + c^2 (N_0 + N'_0 + N_{1/2} + N'_{1/2} )\\
\hskip 1.5 truecm  + r_1 + o(1) \Big ) \leq N_4\\ \\ (N'_5
\leq ) \ 4 \Big ( N'_3 + r_1 + o(1)\Big ) \leq N_5 \end{array}
\right . \eeq

\noi where we have omitted the dependence of the $o(1)$ terms on $N_i$ and $N'_i$. We know in addition that
\beq
\label{2.70e}
N'_{1/2} \vee \widetilde{N}'_{1/2} \leq 2 \Big ( N'_0 (N'_3 + r_1 + o(1))\Big )^{1/2}\ .
\eeq

\noi In order to ensure (\ref{2.69e}), we proceed as follows. We first choose $N_0$ and $N_2$ by imposing
\beq
\label{2.71e}
\left \{ \begin{array}{l} 
N_0 = C_0 \left ( c_4\ N_2 + c_3^2\ N_0 + r_1 + 1 \right )\\ \\ N_2 = C_2 \left ( c_4 \ N_0 + r_2 + 1 \right )
\end{array} \right .
\eeq

\noi which is possible under the smallness condition on $c_3$ , $c_4$
\beq
\label{2.72e}
C_0 \left ( C_2 \ c_4^2 + c_3^2 \right ) < 1\ ,
\eeq

\noi and we impose the condition that $o(1) \leq 1$ in (\ref{2.43e})
(\ref{2.47e}) by taking $T$ sufficiently large (depending on $N_0$,
$N_2$ just chosen and on $N_{1/2}$, $N_1$ to be chosen later). This
ensures the $N'_0 \leq N_0$ and $N'_2 \leq N_2$ parts of (\ref{2.69e}).
Furthermore one can replace $N'_0$ by $N_0$ in all the remaining
estimates. We next impose
\beq
\label{2.73e}
\left \{ \begin{array}{l} 
N_5 = 4 \left ( N_3 + r_1 + 1 \right )\\ \\ N_6 = C_6 \left ( c_4 (N_0\ N_5)^{1/2} + r_2 + 1 \right ) = C_6 \left ( 2c_4 \left ( N_0 (N_3 + r_1 + 1 )\right )^{1/2} + r_2 + 1 \right )
\end{array} \right .
\eeq 

\noi and we impose $o(1) \leq 1$ in (\ref{2.50e}) and (\ref{2.60e}) by
taking $T$ sufficiently large depending on the relevant $N_i$. This
ensures the $N'_6 \leq N_6$ part of (\ref{2.69e}) together with the
inequality
\beq
\label{2.74e}
N'_5 \leq 4 \left ( N'_3 + r_1 + 1 \right )
\eeq 

\noi which will ensure the $N'_5 \leq N_5$ part of (\ref{2.69e}) as
soon as the $N'_3 \leq N_3$ part holds. Furthermore, under the choices
and assumptions made so far, (\ref{2.74e}) implies
\beq
\label{2.75e}
N'_{1/2} \vee \widetilde{N}'_{1/2} \leq 2 \left ( N_0 (N'_3 + r_1 + 1)\right )^{1/2} \ .
\eeq

\noi We now substitute (\ref{2.75e}) into (\ref{2.44e}), we substitute
(\ref{2.74e}) (\ref{2.75e}) into (\ref{2.62e}) and we substitute the
results and (\ref{2.75e}) again into (\ref{2.54e}), thereby obtaining an
estimate of the form 
\bea
\label{2.76e}
N'_3 &\leq& f\left ( N'_3, \{ N_i\}\right ) \equiv C_3 \Big ( 2a \left ( N_0 (N'_3 + r_1 + 1)\right )^{1/2}\nn \\
&&+c_4 \ C_6 \left ( 2c_4 \left (N_0 (N_3 + r_1 + 1)\right )^{1/2} + r_2 + 1\right ) \nn \\
&&+c\ N_2 + c_3^2\ N_3 + 2c^2\ N_0 + r_1 + o(1) \Big ) 
\eea 

\noi where $f(N'_3, \{N_i\})$ is a positive increasing concave function of $N'_3$ for fixed $T$ and $N_i$. It follows therefrom that (\ref{2.76e}) will imply $N'_3 \leq N_3$ provided we ensure that
\beq
\label{2.77e}
N_3 \geq f\left ( N_3, \{N_i\}\right ) \ .
\eeq

\noi This is obtained by imposing 
\bea
\label{2.78e}
N_3 &=& C_3 \Big ( 2\left ( a + C_6 \ c_4^2 \right ) \left ( N_0 (N_3 + r_1 + 1)\right )^{1/2} + c_4\ C_6 (r_2 + 1)\nn \\
&&+ c\ N_2 + c_3^2\ N_3 + 2c^2 \ N_0 + r_1 + 1 \Big )
\eea

\noi which is possible under the smallness condition $C_3c_3^2 < 1$,
and by imposing that $o(1) \leq 1$ in (\ref{2.76e}) by taking $T$
sufficiently large depending on the $N_i$.\par

It is then a simple matter to choose $N_1$ and $N_4$ in order to ensure
the $N'_1 \leq N_1$ and $N'_4 \leq N_4$ parts of (\ref{2.69e}), since
all the $N'_i$ in the RHS of (\ref{2.44e}) and (\ref{2.62e}) are now
under control. It suffices to choose 
\beq
\label{2.79e}
N_1 = C_1 \left ( c_4\ N_2 + 2c^2\ N_0 + a(N_0\ N_5)^{1/2} + r_1 + 1 \right )
\eeq
\beq
\label{2.80e}
N_4 = C_4 \left ( a\ N_5 + (a + 2c^2)(N_0\ N_5)^{1/2} + c(N_6 + N_2) + 2c^2 N_0 + r_1 + 1 \right )
\eeq

\noi and to impose that $o(1) \leq 1$ in (\ref{2.45e}) (\ref{2.63e}) by
taking $T$ sufficiently large depending on the $N_i$ (with the $N'_i$
in the $o(1)$ terms being estimated by the $N_i$).\par

We now show that the map $\phi$ is a contraction on ${\cal R}$ for a
suitable norm defined on $X(I)$. Let $(v_i, B_i) \in {\cal R}$ and
$(v'_i, B'_i) = \phi ((v_i, B_i))$, $i = 1,2$. For any pair of
functions $f_1$, $f_2$, we define $f_{\pm} = (1/2) (f_1 \pm f_2)$ so
that $f_1 = f_+ + f_-$, $f_2 = f_+ - f_-$ and $(fg)_{\pm} = f_+ g_{\pm}
+ f_- g_{\mp}$. In particular $u_+ = u_a + v_+$, $u_- = v_-$, $A_+ =
A_a + B_+$, $A_- = B_-$, $(\Delta_A)_- = - 2iB_- \cdot \nabla_{A_+}$,
and similarly for the primed quantities. Since ${\cal R}$ is convex
and stable under $\phi$, $(v_+, B_+)$ and $(v'_+, B'_+)$ belong to
${\cal R}$, namely satisfy (\ref{2.25e})-(\ref{2.31e}). Corresponding
to (\ref{1.7e}), $(v'_-, B'_-)$ satisfies the system
\beq \label{2.81e} \left \{ \begin{array}{l} i \partial_t v'_- = -
(1/2) (\Delta_A)_+ v'_- + g(|u|_+^2)v'_- + i B_- \cdot \nabla_{A_+} (u_a
+ v'_+) \\ \\
\hskip 1.5 truecm + g\left ( 2 {\rm Re} (\overline{u}_a + \overline{v}_+)v_-\right ) (u_a
+ v'_+) \\ \\ \sq B'_- = 2 P \ {\rm Im} \left ( \overline{v}_-
\nabla_{A_+}(u_a + v_+)\right ) - P\ B_- \left ( |u_a + v_+|^2+|v_-|^2 \right )
\ . \end{array} \right . \eeq

Here however, in contrast with the case of the (WS)$_3$ system where
the corresponding map $\phi$ can be shown to be a contraction for the
whole norm of $X(I)$, we encounter a difficulty due to the derivative
coupling in the covariant Laplacian. In fact if $D$ is a differential
operator of order $m$, a straightforward energy estimate of $\parallel
Dv'_-\parallel_2$ from (\ref{2.81e}) yields
$$\partial_t \parallel Dv'_-\parallel_2\ \leq \ \parallel B_-\cdot D\nabla_{A_+} v'_+\parallel_2 + \ \hbox{other terms}$$

\noi and requires therefore a control of $v'_+$ at order $m+1$, so that
one can hope to contract norms of $v$ of degree at most one less than
those occurring in the definition of $X(I)$. Fortunately, because of
the special algebraic properties of the equations, it turns out that
the lowest two semi norms of $X(I)$ for the differences, namely those
corresponding to $N_0$ and $N_2$, can be decoupled from the higher ones
and can be contracted on the bounded sets of $X(I)$. This follows from the fact that
the symmetry of the quadratic form $P\ {\rm Im}(\overline{v}_1
\nabla_A v_2)$ has made it possible to avoid having a gradient acting
on $v_-$ in the equation for $B'_-$ in (\ref{2.81e}). Thus we shall
show that $\phi$ is a contraction for the pair of semi norms
\beq \label{2.82e} \left \{ \begin{array}{l} N_0 = \ \displaystyle{\mathrel{\mathop {\rm Sup}_{t \in I }}}\ h(t)^{-1} \parallel v(t) \parallel_2 \ , \\ \\
N_2 = \ \displaystyle{\mathrel{\mathop {\rm Sup}_{t \in I }}}\ h(t)^{-1} \parallel B;L^4([t, \infty ),L^4 )\parallel \ . \end{array} \right . \eeq

\noi Let $(N_{0-},N_{2-})$ and $(N'_{0-}, N'_{2-})$ be the
corresponding semi norms of $(v_-, B_-)$ and $(v'_-, B'_-)$
respectively. We have to estimate $(N'_{0-}, N'_{2-})$ in terms of
$(N_{0-}, N_{2-})$. We first estimate $N'_{0-}$. From (\ref{2.81e}) we
obtain 
\beq
\label{2.83e}
\parallel v'_-(t)\parallel_2\ \leq\ \parallel B_- \cdot \nabla_{A_+}(u_a
+ v'_+)\parallel_+ \ + \ \parallel g\left ( 2 {\rm Re} \left ( \overline{u}_a
+ \overline{v}_+\right ) v_- \right ) (u_a + v'_+)\parallel_+ \ . \eeq

\noi The terms not containing $v'_+$ are estimated as in the proof of (\ref{2.42e}), namely 
\begin{eqnarray*}
&&\parallel B_- \cdot \nabla u_a \parallel_+ \ \leq C \ c_4 \ N_{2-}\ h(t) \ , \\
&&\parallel B_- \cdot A_a u_a \parallel_+ \ \leq c\ a \ N_{2-}\ t^{-1} \ h(t) \ ,\\
&&\parallel B_- B_+ u_a \parallel_+ \ \leq c\ N_{2-} \ N_2 \ t^{-1} \ h(t)^2 \ , 
\end{eqnarray*}
$$\parallel g\left ( \overline{u}_a v_-\right ) u_a \parallel_+ \ \leq C\ c_3^2\ N_{0-}\ h(t) \ ,$$
$$\parallel g\left ( \overline{v}_+ v_-\right ) u_a\parallel_+ \ \leq C\ c_3\ N_{0-}\ N_1\ t^{1/4}\ h(t)^2 \ .$$

\noi We next estimate the terms containing $v'_+$.
$$\parallel B_- \cdot \nabla v'_+ \parallel_+ \ \leq\ \parallel \ \parallel B_- \parallel_4\ \parallel \nabla v'_+ \parallel_4 \ \parallel_1 \ \leq C\ N_{2-}\ N_4\ h(t) \ \overline{h}(t)$$

\noi by Lemma 2.2,
$$\parallel B_- \cdot A_a v'_+ \parallel_+ \ \leq\ \parallel\ \parallel B_- \parallel_4\ \parallel A_a \parallel_{\infty} \ \parallel v'_+ \parallel_4 \ \parallel_1$$
$$\leq a\ N_{2-}\ N_1\ h^2 \parallel t^{-1} \parallel_{8/3}\ \leq a \ N_{2-} \ N_1\ t^{-5/8} \ h(t)^2 \ ,$$
$$\parallel B_- B_+ v'_+ \parallel_+ \ \leq \ C \parallel \ \parallel B_-\parallel_4 \left ( \parallel B_+ \parallel_4 \ \parallel \nabla B_+ \parallel_4^3 \right )^{1/4} \parallel v'_+ \parallel_4 \ \parallel_1$$
$$\leq C\ N_{2-}(N_2\ N_6^3)^{1/4} \ N_1 \ t^{1/8}\ h(t)^3$$

\noi by Lemma 2.2, 
\begin{eqnarray*}
\parallel g\left ( \overline{u}_a v_- \right ) v'_+ \parallel_+ &\leq&  C \parallel\ \parallel u_a \parallel_3\ \parallel v_- \parallel_2 \ \parallel v'_+ \parallel_3\ \parallel_1 \\
&\leq& C\ c_3\ N_{0-}\ N_1 \ t^{1/4}\ h(t)^2 \ ,
\end{eqnarray*}
\begin{eqnarray*}
\parallel g\left ( \overline{v}_+ v_- \right ) v'_+ \parallel_+ &\leq&  C \parallel\ \parallel v_+ \parallel_3\ \parallel v_- \parallel_2 \ \parallel v'_+ \parallel_3\ \parallel_1 \\
&\leq& C\ N_{0-}\ N_1^2 \ t^{1/2}\ h(t)^3 \ ,
\end{eqnarray*}
 
\noi by Lemma 2.3, part (1) and Lemma 2.2.\par

Collecting the previous estimates yields
\bea
\label{2.84e}
&&N'_{0-} \leq C_0 \Big \{ N_{2-} \Big ( c_4 + a\ c\ T^{-1} + c\ N_2\ T^{-1}\ h(T) + N_4 \ \overline{h}(T) + a\ N_1\ T^{-5/8} \ h(T)\nn \\
&&+ (N_2\ N_6^3)^{1/4} \ N_1\ T^{-1/4} \ h(T) \ \overline{h}(T)\Big ) + N_{0-} \Big ( c_3^2 + c\ N_1 \ T^{-1/8} \ \overline{h}(T)\nn \\
&&+ N_1^2\ T^{-1/4}\ \overline{h} (T)^2 \Big ) \Big \}
\eea

\noi which is of the form 
\beq
\label{2.85e}
N'_{0-} \ \leq C_0 \Big ( N_{2-} (c_4 + o(1) ) + N_{0-} ( c_3^2 + o(1))\Big ) \ .
\eeq

\noi We next estimate $N'_{2-}$. From (\ref{2.8e}) (\ref{2.81e}) we obtain
\beq
\label{2.86e}
\parallel B'_-;L^4(J,L^4)\parallel\ \leq \ C \left ( \parallel \overline{v}_- \nabla_{A_+} (u_a + v_+) \parallel_* \ + \ \parallel B_- (|u_a|^2 + |v_+|^2 )\parallel_* \right )\ .
\eeq

\noi The linear terms are estimated as in the proof of (\ref{2.46e}), namely 
\begin{eqnarray*}
&&\parallel \overline{v}_- \nabla u_a \parallel_*\ \leq 2c_4\ N_{0-}\ h(t) \ ,\\
&&\parallel \overline{v}_- A_a u_a \parallel_*\ \leq a\ c\ N_{0-}\ t^{-1}\ h(t) \ ,\\
&&\parallel B_- |u_a|^2 \parallel_*\ \leq c^2\ N_{2-}\ t^{-1}\ h(t) \ .\\
\end{eqnarray*}
\noi The non linear terms are estimated in a slightly different way. The quadratic terms are estimated by
$$\parallel \overline{v}_- \nabla v_+ \parallel_* \ \leq \ \parallel\ \parallel v_- \parallel_2\ \parallel \nabla v_+ \parallel_4\ \parallel_{4/3} \ \leq C\ N_{0-}\ N_4\ h(t) \ \overline{h}(t) \ ,$$
$$\parallel \overline{v}_- A_a v_+ \parallel_* \ \leq \ \parallel\ \parallel A_a \parallel_{\infty}\ \parallel v_- \parallel_2\ \parallel v_+ \parallel_4 \ \parallel_{4/3} \ \leq a\ N_{0-}\ N_1\ t^{-5/8}\ h(t)^2  \ ,$$
$$\parallel \overline{v}_- B_+ u_a\parallel_* \ \leq \ \parallel\ \parallel v_- \parallel_2\ \parallel B_+ \parallel_4\ \parallel u_a \parallel_{\infty} \ \parallel_{4/3}$$
$$\leq c\ N_{0-}\ N_2\ h \parallel t^{-3/2} h \parallel_2\ \leq c\ N_{0-}\ N_2\ t^{-1} \ h(t)^2 \ . $$

\noi The cubic terms are estimated by
$$\parallel \overline{v}_- B_+ v_+\parallel_* \ \leq \ C \parallel\ \parallel v_- \parallel_2\left (  \parallel B_+ \parallel_4\ \parallel \nabla B_+ \parallel_4^3\right )^{1/4} \parallel v_+ \parallel_4\ \parallel_{4/3}$$
$$\leq C\ N_{0-} (N_2 N_6^3)^{1/4} N_1\ t^{1/8} \ h(t)^3$$

\noi by Lemma 2.2,
$$\parallel B_- |v_+|^2\parallel_* \ \leq \ C\ N_{2-}\ N_1^{3/2}\ N_{1/2}^{1/2} \ h(t)^3 \ .$$

\noi Collecting the previous estimates yields
$$N'_{2-} \leq C_2 \Big \{ N_{0-} \left (  c_4 + a\ c \ T^{-1} + N_4\ \overline{h}(T) + a\ N_1\ T^{-5/8} \ h(T) + c \ N_{2} \ T^{-1} h(T) \right .$$
\beq
\label{2.87e}
\left. + (N_2\ N_6^3)^{1/4} \ N_1\  T^{-1/4} \ h(T) \ \overline{h}(T) \right ) 
+ N_{2-} \left ( c^2\ T^{-1} + N_1^{3/2}\ N_{1/2}^{1/2}\ h(T)^2 \right ) \Big \}
\eeq

\noi which is of the form 
\beq
\label{2.88e}
N'_{2-} \leq C_2\left ( c_4 + o(1) \right )  N_{0-} + o(1)\ N_{2-}\ .
\eeq

\noi It follows from (\ref{2.85e}) (\ref{2.88e}) that the map $\phi$ is a contraction for the pair of semi norms $(N_0, N_2)$ on the set ${\cal R}$ under the smallness condition
\beq
\label{2.89e}
C_0 \left ( C_2\ c_4^2 + c_3^2 \right ) < 1
\eeq

\noi and for $T$ sufficiently large. Since the set ${\cal R}$ is closed
for the norm defined by the pair $(N_0, N_2)$, it follows therefrom
that the system (\ref{1.4e}) has a solution in ${\cal R}$. This proves
the existence part of Proposition 2.2. The uniqueness part follows from
(\ref{2.85e}) (\ref{2.88e}) again with $N'_{i-} = N_{i-}$. \par

We remark at this point that the constants $C_0$, $C_2$ appearing in
(\ref{2.89e}) can be taken to be the same as in (\ref{2.72e}) so that
the two smallness conditions actually coincide. In fact those constants
are determined by the linear terms in the estimates, and those terms
are the same in both cases. There may occur additional, different
constants coming from the non linear terms. They have been omitted in
(\ref{2.84e}) (\ref{2.87e}). \par

It remains to prove the last statement of Proposition 2.2 and for that purpose we need to estimate the energy norm of $B'$. From (\ref{1.7e}) (\ref{2.10e}) it follows that for all $t \in I$
\beq
\label{2.90e}
\parallel \nabla B'(t) \parallel_2 \ \vee ÷ \parallel \partial_t B'(t) \parallel_2 \ \leq \ \parallel G_2 - R_2 \parallel_+
\eeq

\noi where $G_2$ is defined by (\ref{1.5e}). We estimate the various terms of $G_2$ successively. The linear terms in $v$ are estimated by 
\begin{eqnarray*}
\parallel \overline{v} \nabla u_a\parallel_+ &\leq& \parallel \ \parallel v \parallel_2\  \parallel \nabla u_a \parallel_{\infty} \ \parallel_{1} \\
&\leq& c\ N_0 \parallel t^{-3/2} h \parallel_{1} \ \leq \ 2c\ N_0\ t^{-1/2}\ h(t) \ ,
\end{eqnarray*}
\begin{eqnarray*}
\parallel \overline{v} A_a u_a\parallel_+ &\leq& \parallel \ \parallel v \parallel_2\  \parallel A_a \parallel_{\infty}\  \parallel u_a \parallel_{\infty}\ \parallel_{1} \\
&\leq& a\ c\ \ N_0 \parallel t^{-5/2} h \parallel_{1} \ \leq \ a\ c\ N_0\ t^{-3/2}\ h(t) \ .
\end{eqnarray*}

\noi The linear term in $B$ is estimated by
\begin{eqnarray*}
\parallel B| u_a|^2 \parallel_+ &\leq& \parallel \ \parallel B \parallel_4\  \parallel u_a \parallel_4\ \parallel u_a \parallel_{\infty}\ \parallel_{1} \\
&\leq& c^2\ N_2\ h  \parallel t^{-3/4-3/2}\parallel_{4/3} \ \leq \ c^2\ N_2\ t^{-3/2}\ h(t) \ .
\end{eqnarray*}

\noi The quadratic terms in $v^2$ are estimated by
$$\parallel \overline{v} \nabla v \parallel_+ \ \leq \ \parallel \ \parallel v \parallel_4\  \parallel \nabla v \parallel_{4}\   \parallel_{1} \ \leq \  C\  N_1\ \ N_{4} \ t^{1/4}  \ h(t)^2$$

\noi by Lemma 2.2,
\begin{eqnarray*}
\parallel A_a| v|^2 \parallel_+ &\leq& \parallel \ \parallel A_a \parallel_{\infty}\  \parallel v \parallel_4^2\  \parallel_1 \\
&\leq& a\ N_1^2\ h^2  \parallel t^{-1}\parallel_{4} \ \leq \ a\ N_1^2\ t^{-3/4}\ h(t)^2 \ .
\end{eqnarray*}

\noi The quadratic terms in $Bv$ again need not be considered. The cubic term $B|v|^2$ is estimated by
$$\parallel B|v|^2\parallel_+\   \leq \ C \parallel B;L^4(L^4) \parallel\ \parallel v;L^{8/3}(L^4)\parallel^{5/4} \ \parallel \nabla v;L^{8/3}(L^4)\parallel^{3/4}$$
$$\leq C\ N_2\ N_1^{5/4}\ N_4^{3/4}\ h(t)^3\ .$$

Collecting the previous estimates and using (\ref{2.23e}), we obtain
\bea
\label{2.91e}
&&\parallel \nabla B'(t)\parallel_2\ \vee \ \parallel \partial_t B'(t) \parallel_2\ \leq C\Big ( c\ N_0\ t^{-1/2} + a\ c\ N_0\ t^{-3/2} + c^2\ N_2 \ t^{-3/2}\nn \\
&&+ N_1\ N_4\ t^{1/4}\ h(t) + a\ N_1^2\ t^{-3/4} \ h(t) + N_2\ N_1^{5/4}\ N_4^{3/4}\ h(t)^2 + r_2 \ t^{-1/2}\Big ) h(t)\nn \\
\eea

\noi which proves that the solution of (\ref{1.4e}) constructed previously satisfies (\ref{2.24e}). \par \nobreak \hfill $\sq$ \par

\noi {\bf Remark 2.2.} The only smallness condition on $u$ is the
condition (\ref{2.72e}), coming from $N_0$ and from its coupling with $N_2$.
The subsequent condition $C_3c_3^2 < 1$ needed for the choice of $N_3$
comes in fact from exactly the same estimate as the $c_3^2$
contribution to $N'_0$, so that the latter condition is actually the
$c_4 = 0$ special case of (\ref{2.72e}) and is therefore weaker than
(\ref{2.72e}). That fact is hidden by the use of overall constants
$C_0$ and $C_3$ in the estimates of $N'_0$ and $N'_3$.

\mysection{Remainder estimates and completion of the proof}
\hspace*{\parindent}
In this section, we first prove that the choice of asymptotic functions
$(u_a,A_a)$ made in the introduction satisfies the assumptions of
Proposition 2.2 for the choice of $h$ made in Proposition 1.1, under
suitable assumptions on the asymptotic state $(u_+, A_+, \dot{A}_+)$.
We then combine those results with Proposition 2.2 to complete the
proof of Proposition 1.1. \par

We first supplement the definition of $(u_a, A_a)$ with some additional
properties of a general character. In addition to the representation
(\ref{1.13e}) (\ref{1.14e}) of $A_1$, we need a representation of
$\partial_t A_1$. From (\ref{1.12e}) it follows that 
\beq
\label{3.1e}
\partial_t A_1(t) = - \int_t^{\infty} dt'\ \cos (\omega (t' - t)) t'^{-1} \ P\ x |u_a(t')|^2
\eeq

\noi so that upon substitution of (\ref{1.8e}) we obtain
\beq
\label{3.2e}
\partial_t A_1(t) = t^{-2} \ D_0 (t) \ \widetilde{\widetilde{A}}_1
\eeq

\noi where
\beq
\label{3.3e}
\widetilde{\widetilde{A}}_1 = - \int_t^{\infty} d\nu \ \nu^{-3} \cos (\omega (\nu - 1)) D_0(\nu ) \ P \ x |w_+|^2 \ .
\eeq

\noi On the other hand, from (\ref{1.13e})
\beq
\label{3.4e}
\nabla A_1(t) = t^{-2} D_0 (t) \nabla \widetilde{A}_1 \ .
\eeq

\noi We shall need the operator 
\beq
\label{3.5e}
J \equiv J(t) = x + it \nabla \ .
\eeq

\noi The asymptotic form $A_a$ for $A$ has been chosen in order to make $R_2$ small. In fact $R_2$ can be rewritten as
\beq
\label{3.6e}
R_2 = \sq A_a + P \left ( t^{-1} \ {\rm Re}\ \overline{u}_a \ J u_a + (A_a - x/t) |u_a|^2 \right )
\eeq

\noi and $A_a$ has been chosen in such a way that 
\beq
\label{3.7e}
\sq A_a = P(x/t)|u_a|^2
\eeq

\noi so that
\beq
\label{3.8e}
R_2 = P \left ( t^{-1} \ {\rm Re}\ \overline{u}_a\ J u_a + A_a |u_a|^2 \right ) \ .
\eeq

Under general assumptions on $(u_a, A_a)$, of the same type as in
Proposition 2.2 (see especially (\ref{2.17e}) (\ref{2.20e})) but not
making use of their special form, we can prove that $R_2$ satisfies
the assumptions needed for that proposition with the choice of $h$
required for Proposition 1.1. \\

\noi {\bf Lemma 3.1.} {\it Let $(u_a, A_a)$ satisfy the estimates
\beq
\label{3.9e}
\parallel u_a(t)\parallel_r \ \leq c\ t^{-\delta (r)} \qquad \hbox{\it for $2 \leq r \leq \infty$} \ ,
\eeq
\beq
\label{3.10e}
\parallel \nabla u_a(t)\parallel_4 \ \leq c\ t^{-3/4} \ ,
\eeq
\beq
\label{3.11e}
\parallel J u_a(t)\parallel_2 \ \leq c_1\ (1 + \ell n\ t) \ ,
\eeq
\beq
\label{3.12e}
\parallel A_a(t)\parallel_{\infty} \ \vee \ \parallel \nabla A_a(t)\parallel_{\infty}\ \leq a\ t^{-1} 
\eeq

\noi for all $t \geq 1$. Then $R_2$ satisfies the estimates
\beq
\label{3.13e}
\parallel R_2;L^{4/3}(J,L^{4/3})\parallel \ \vee \ \parallel \nabla R_2;L^{4/3}(J,L^{4/3})\parallel \ \leq r_2\ t^{-1}(1 + \ell n \ t) \ , 
\eeq
\beq
\label{3.14e}
\parallel R_2;L^1(J,L^2)\parallel \ \leq r_2\ t^{-3/2} (1 + \ell n \ t) 
\eeq

\noi for some constant $r_2$ and for all $t \geq 1$, where $J = [t, \infty )$.}\\

\noi {\bf Proof.} We estimate
\begin{eqnarray*}
\parallel R_2(t)\parallel_{4/3} &\leq& C\parallel u_a\parallel_4 \left ( t^{-1} \parallel Ju_a \parallel_2 \ + \ \parallel A_a \parallel_{\infty} \ \parallel u_a \parallel_2\right )\\
&\leq& C\ t^{-7/4}\ c \left ( c_1 (1 + \ell n\ t) + ac\right )
\end{eqnarray*}

\noi which implies the first estimate of (\ref{3.13e}) by integration, 
\begin{eqnarray*}
\parallel R_2(t)\parallel_2 &\leq& \parallel u_a\parallel_{\infty} \left ( t^{-1} \parallel Ju_a \parallel_2 \ + \ \parallel A_a \parallel_{\infty} \ \parallel u_a \parallel_2\right )\\
&\leq& t^{-5/2}\ c \left ( c_1 (1 + \ell n\ t) + ac\right )
\end{eqnarray*}

\noi which implies (\ref{3.14e}) by integration.\par

In order to prove the second estimate of (\ref{3.13e}), we note that the quadratic form 
$$P\ {\rm Re} \left ( t^{-1} \ \overline{v}_1 \ J v_2 + A_a \ \overline{v}_1 v_2 \right )$$

\noi is symmetric in $v_1$, $v_2$, so that 
\begin{eqnarray*}
\nabla  R_2 &=& 2\ P\ {\rm Re} \left ( t^{-1} (\nabla  \overline{u}_a) \ Ju_a + A_a (\nabla \overline{u}_a) u_a \right )\\
&& + P \left ( \nabla (A_a + x/t)\right ) |u_a|^2
\end{eqnarray*}

\noi and therefore
\begin{eqnarray*}
\parallel \nabla R_2\parallel_{4/3} &\leq& C\Big ( \parallel \nabla u_a\parallel_{4} \left ( t^{-1} \parallel Ju_a \parallel_2 \ + \ \parallel A_a \parallel_{\infty} \ \parallel u_a \parallel_2\right )\\
&&+ \left ( \parallel \nabla A_a \parallel_{\infty}  + t^{-1}\right ) \parallel u_a \parallel_4\ \parallel u_a \parallel_2 \Big ) \\
&\leq& C\ t^{-7/4}\left ( c \left ( c_1 (1 + \ell n\ t) + ac\right ) + c^2 (a+1) \right )
\end{eqnarray*}

\noi from which the second estimate of (\ref{3.13e}) follows by integration. \par \nobreak \hfill $\sq$ \par

We now turn to $R_1$. We first skim $R_1$ of some harmless terms.
Expanding the covariant Laplacian and using again $J$, we rewrite
$R_1$ as
\beq
\label{3.15e}
R_1 = R_{1,1} + R_{1,2}
\eeq

\noi where
\beq
\label{3.16e}
R_{1,1} = i \partial_t u_a + (1/2) \Delta u_a + t^{-1} (x \cdot A_1) u_a - g(|u_a|^2)u_a \ ,
\eeq
\beq
\label{3.17e}
R_{1,2} = t^{-1}(x \cdot A_0) u_a - t^{-1} A_a \cdot Ju_a - (1/2) A_a^2 u_a \ .
\eeq

\noi In the same way as for $R_2$, we can show that $R_{1,2}$ satisfies
the assumptions needed for Proposition 2.2 with the choice of $h$
required for Proposition 1.1 under general assumptions on $(u_a, A_a)$
not making use of their special form.\\

\noi {\bf Lemma 3.2.} {\it Let $u_a$, $A_a$ and $A_0$ satisfy the estimates
\beq
\label{3.18e}
\parallel \partial_t^j \nabla^k u_a \parallel_2\ \leq c \ ,
\eeq
\beq
\label{3.19e}
\parallel \partial_t^j \nabla^k J u_a \parallel_2\ \leq c_1(1 + \ell n\ t) \ ,
\eeq
$$\parallel \partial_t^j \nabla^k A_a \parallel_{\infty}\ \leq a\ t^{-1} \ , \eqno(2.20)\equiv(3.20)$$
$$\parallel \partial_t^j \nabla^k (x \cdot A_0) \parallel_{\infty}\ \leq a_0\ t^{-1} \ , \eqno(3.21)$$

\noi for $0 \leq j + k \leq 1$ and for all $t \geq 1$. Then $R_{1,2}$ satisfies the estimates 
$$\parallel \partial_t^j \nabla^k R_{1,2} \parallel_{2}\ \leq r_{1,2} \ t^{-2} (1 + \ell n \ t)\ , \eqno(3.22)$$

\noi for $0 \leq j + k \leq 1$, for some constant $r_{1,2}$ and for all $t \geq 1$.}\\

\noi {\bf Proof.} We estimate
$$\parallel R_{1,2} \parallel_2 \ \leq t^{-2} \left ( (a_0 + (1/2) a^2)c + ac_1(1 + \ell n\ t)\right )\ ,$$
$$\parallel \nabla R_{1,2} \parallel_2 \ \leq t^{-2} \left ( (2a_0 + (3/2) a^2)c + 2ac_1(1 + \ell n\ t)\right )\ ,$$
$$\parallel \partial_t R_{1,2} \parallel_2 \ \leq \hbox{\rm idem} + t^{-3} \left ( (a_0  + ac_1(1 + \ell n\ t)\right )\ .$$
\par \nobreak \hfill $\sq$\par

We now turn to $R_{1,1}$. We shall need the commutation relations 
$$\nabla MD = MD \left ( ix + t^{-1} \nabla \right ) \equiv MD \widetilde{\nabla} \ , \eqno(3.23)$$
$$i \partial_t MD = MD \left ( i \partial_t + (1/2) x^2 - i t^{-1} (x \cdot \nabla + 3/2)\right ) \equiv MD\ i \widetilde{\partial}_t \ , \eqno(3.24)$$
$$JMD = i\ MD \ \nabla \ , \eqno(3.25)$$
$$\left ( i \partial_t + (1/2) \Delta \right ) MD = MD \left ( i \partial_t + (2t^2)^{-1} \Delta \right )\ . \eqno(3.26)$$

In particular (3.23) (3.24) are taken as the definitions of
$\widetilde{\nabla}$ and $\widetilde{\partial}_t$. From the choice
(\ref{1.8e}) of $u_a$ and from (3.26), it follows that 
$$R_{1,1} = MD \left ( i \partial_t + (2t^2)^{-1} \Delta + t^{-1} x \cdot \widetilde{A}_1 - t^{-1} g(|w_+|^2)\right ) \exp (-i \varphi ) w_+ \ . \eqno(3.27)$$

\noi The choice (\ref{1.15e}) of $\varphi$ has been taylored to cancel the two long range terms in (3.27), so that
$$R_{1,1} = (2t^2)^{-1} \ MD\ \Delta \exp (- i \varphi ) w_+ \ . \eqno(3.28)$$

We now have to prove that the previous choice of $(u_a, A_a)$ satisfies
the remaining assumptions of Proposition 2.2 and of Lemmas 3.1 and 3.2.
More precisely we have to prove that $(u_a, A_a)$ satisfies the
estimates (\ref{2.17e}) (\ref{2.19e}) (\ref{2.20e}) (\ref{3.19e})
(3.21) and the analogue of (3.22) for $R_{1,1}$. (Note that
(\ref{3.9e}) (\ref{3.10e}) (\ref{3.18e}) are special cases of
(\ref{2.17e}) and that (\ref{3.12e}) is a special case of (3.20) which
is identical with (\ref{2.20e})).\par

The contribution of $A_0$ to $A_a$ and to $R_{1,2}$ will be taken care
of by the following general estimates of solutions of the wave
equation. \\

\noi{\bf Lemma 3.3.} {\it Let $A_0$ be defined by (\ref{1.11e}) and let $k \geq 0$ be an integer. Let $A_+$ and $\dot{A}_+$ satisfy the conditions
$$\nabla^2A_+\ , \ \nabla \dot{A}_+ \in W_1^k \qquad , \quad A_+ \in L^3 \ , \ \dot{A}_+ \in L^{3/2} \ . \eqno(3.29)$$

\noi Then $A_0$ satisfies the estimates}
$$\left \{ \begin{array}{l} \parallel A_0(t); W_{\infty}^k\parallel \ \leq \ a_0 \ t^{-1}\\ \\ \parallel \partial_t A_0(t); W_{\infty}^{k-1}\parallel \ \leq \ a_0 \ t^{-1}\qquad \hbox{\it for $k \geq 1$} \ .\end{array} \right . \eqno(3.30)$$

A proof can be found in \cite{20r}. As mentioned in Remark 1.2, the
assumptions $A_+ \in L^3$ and $\dot{A}_+ \in L^{3/2}$ serve to
exclude constants in $A_+$ and $\dot{A}_+$ and linear terms in $x$ in
$A_+$, but are otherwise controlled by the $W_1^k$ assumption through
Sobolev inequalities.\\

We next derive some preliminary estimates of $\widetilde{A}_1$ and $\widetilde{\widetilde{A}}_1$.\\

\noi{\bf Lemma 3.4.} {\it Let $k \geq 0$ be an integer. Then the following estimates hold.} 
$$\parallel \omega^{k+1}  \widetilde{A}_1 \parallel_2\ \vee \ \parallel \omega^k \widetilde{\widetilde{A}}_1\parallel_2 \ \leq (k+1/2)^{-1} \ \parallel \omega^k x|w_+|^2 \parallel_2 \ , \eqno(3.31)$$
$$\parallel \omega^{k+1}(x\cdot   \widetilde{A}_1) \parallel_2\ \leq \ (k-1/2)^{-1}\left ( \parallel \omega^k x^2|w_+|^2\parallel_2 \ + \ 2\parallel \omega^k x|w_+|^2 \parallel_2\right ) \ {\it for}\ k\geq 1 \ , \eqno(3.32)$$
$$\parallel \nabla^k \widetilde{A}_1 \parallel_{\infty}\ \leq \ C\parallel \omega^k x|w_+|^2;H^1\parallel \ , \eqno(3.33)$$
$$\parallel \nabla^{k-1} \widetilde{\widetilde{A}}_1 \parallel_{\infty}\ \leq \ C\parallel \omega^k x|w_+|^2;H^1\parallel \quad \hbox{\it for $k \geq 1$}\ , \eqno(3.34)$$
$$\parallel \nabla^k (x \cdot \widetilde{A}_1) \parallel_{\infty}\ \leq \ C\left ( \parallel \omega^k x^2|w_+|^2;H^1\parallel \ + \ \parallel \omega^k x|w_+|^2;H^1\parallel\right )\ \hbox{\it for $k \geq 1$}\ .  \eqno(3.35)$$

\noi {\bf Proof.} (3.31) follows immediately from (\ref{1.14e}) and (\ref{3.3e}). From (\ref{1.14e}) and from the commutation relation 
$$[x;P] = - 2 \omega^{-2}\nabla$$

\noi it follows that
$$x\cdot \widetilde{A}_1 = \int_1^{\infty} d\nu \ \nu^{-2} \ \omega^{-1} \sin (\omega (\nu - 1)) D_0(\nu ) \left \{ P \cdot (x \otimes x) |w_+|^2 - 2 \omega^{-2} \nabla \cdot x |w_+|^2 \right \} \eqno(3.36)$$

\noi so that
$$\parallel \omega^{k+1}\ x\cdot   \widetilde{A}_1 \parallel_2\ \leq \ \int_1^{\infty} d\nu \ \nu^{-2} \Big ( \parallel \omega^{k}D_0(\nu) x^2 |w_+|^2 \parallel_2 $$
$$+ 2(\nu - 1) \parallel \omega^{k+1} D_0 (\nu ) \omega^{-1} x|w_+|^2 \parallel_2 \Big )$$
$$\leq \int_1^{\infty} d\nu \ \nu^{-1/2-k} \left ( \parallel \omega^k x^2|w_+|^2 \parallel_2 \ + \ 2 \parallel \omega^k x|w_+|^2 \parallel_2 \right )$$

\noi which implies (3.32). Finally (3.33)-(3.35) follow from (3.31) (3.32) and from the fact that $\dot{H}^1 \cap \dot{H}^2 \subset L^{\infty}$. \par \nobreak \hfill $\sq$\par

As an immediate corollary, we obtain the following estimates of $A_1$ and $\partial_t A_1$.\\

\noi{\bf Corollary 3.1.} {\it The following estimates hold.}
$$\parallel A_1(t) \parallel_{\infty} \ \leq \ C\parallel w_+;H^{2,1} \parallel^2 \ t^{-1} \ , \eqno(3.37)$$
$$\parallel \partial_t A_1(t) \parallel_{\infty} \ \vee \ \parallel \nabla A_1(t) \parallel_{\infty}\ \leq \ C\parallel w_+;H^{2,1} \parallel^2 \ t^{-2} \ . \eqno(3.38)$$

\noi {\bf Proof.} The result follows from (\ref{1.13e}) (\ref{3.2e}) (\ref{3.4e}) and from (3.33) (3.34). \par \nobreak \hfill $\sq$ \par

We next derive the remaining estimates of $u_a$ and of $R_{1,1}$. The
following proposition is slightly stronger than needed.\\

\noi {\bf Proposition 3.1.} {\it Let $u_a$ be defined by (\ref{1.8e})
with $w_+ = Fu_+$ and with $\varphi$ defined by (\ref{1.16e})
(\ref{1.2e}) (\ref{1.14e}) and let $R_{1,1}$ be given by (3.28). Let
$u_+ \in H^{3,1} \cap H^{1,3}$. Then the following estimates hold for
some constants $c$, $c_1$ and $r_{1,1}$, for $0 \leq j+k \leq 1$ and
for all $t\geq 1$~: 
$$\parallel \partial_t^j \nabla^k u_a(t)\parallel_r \ \leq c\ t^{-\delta (r)} \qquad \hbox{\it for $2 \leq r \leq \infty$}\ . \eqno(2.17)\equiv(3.39)$$

\noi In particular}
$$\parallel u_a(t) \parallel_3\ \leq \ \parallel w_+ \parallel_3 \ t^{-1/2} \ , \eqno(3.40)$$
$$\parallel \nabla u_a(t) \parallel_4\ \leq \left (  \parallel xw_+ \parallel_4\ + O(t^{-1} \ell n\ t)\right ) t^{-3/4} \ .   \eqno(3.41)$$
$$\parallel \partial_t^j \nabla^{k+1} u_a(t) \parallel_r\ \leq \ c\ t^{-\delta (r)}\qquad \hbox{\it for $2 \leq r \leq 6$} \ , \eqno(3.42)$$
$$\parallel \partial_t^j \nabla^{k} J u_a(t) \parallel_r\ \leq \ c_1(1 + \ell n\ t)\ t^{-\delta (r)}\qquad \hbox{\it for $2 \leq r \leq 6$} \ , \eqno(3.43)$$
$$\parallel \partial_t^j \nabla^{k} R_{1,1}(t) \parallel_2\ \leq \ r_{1,1}\ t^{-2}(1 + \ell n\ t)^2 \ .  \eqno(3.44)$$

\noi {\bf Proof.} From the commutation relations (3.23) (3.24), it follows that for any differential operator $Z$
$$\parallel \partial_t^j \nabla^{k} MD\ Z \exp (- i \varphi ) w_+ \parallel_r\ = \ t^{-\delta (r)}\parallel \widetilde{\partial}_t^j \ \widetilde{\nabla}^{k}  Z \exp (- i \varphi ) w_+ \parallel_r\ . \eqno(3.45)$$

\noi From (\ref{1.8e}) (3.28) and from the commutation relation (3.25), it follows that in order to derive (3.39)-(3.44) we have to estimate norms of the type 
$$\parallel \widetilde{\partial}_t^j \ \widetilde{\nabla}^{k}  Z \exp (- i \varphi ) w_+ \parallel_r$$

\noi for $0 \leq j + k \leq 1$, for suitable choices of $Z$ and $r$, and with suitable $r$-independent time behaviour. The relevant choices are
\begin{eqnarray*}
&&Z = 1 \ , \ 2 \leq r \leq \infty \quad \hbox{\rm for (3.39)}\\
&&Z = \widetilde{\nabla}\ {\rm or}\ \nabla  \ , \ 2 \leq r \leq 6 \quad \hbox{\rm for (3.42) (3.43)}\ ,\\
&&Z = t^{-2} \Delta  \ , \ r = 2 \quad \hbox{\rm for (3.44)}\ . 
\end{eqnarray*}

\noi Expanding $\widetilde{\partial}_t^j \ \widetilde{\nabla}^k$
according to the definitions (3.23) (3.24) and omitting the commutators of
derivatives with powers of $x$ and $t$ which generate terms of lower
order, we are led to estimate norms of the type $\parallel Z \exp (-i
\varphi ) w_+ \parallel_r$ for the following choices of $Z$, $r$~:
$$Z = 1, x, t^{-1}\nabla , x^2, \partial_t, t^{-1} x \cdot \nabla \hskip 2.3 truecm {\rm with}\ 2 \leq r \leq \infty \ {\rm for} \ (3.39)\ ,$$
$$Z = x, x^2, t^{-1} x \nabla , x^3, x \partial_t, t^{-1} x^2 \nabla, t^{-1} \nabla , t^{-2} \nabla^2, t^{-1} \partial_t \nabla, t^{-2} x \nabla^2$$
\hskip  9 truecm ${\rm with} \ 2 \leq r \leq 6 \ {\rm for} \ (3.42)\ ,$
$$Z = \nabla , x \nabla, t^{-1} \nabla^2, x^2 \nabla, \partial_t \nabla, t^{-1}x\nabla^2\ \hskip 1.2 truecm {\rm with} \ 2 \leq r \leq 6\ {\rm for}\ (3.43)\ , $$
$$Z = \Delta ,x \Delta , t^{-1} \nabla \Delta , x^2\Delta , \partial_t \Delta , t^{-1}x \nabla \Delta\quad  {\rm with} \ r=2\ {\rm for} \ (3.44) \ ,$$

\noi where we have omitted an overall $t^{-2}$ factor in the last case. \par

We expand the derivatives acting on $\exp (- i \varphi ) w_+$ by the
Leibnitz rule and we estimate the expressions thereby obtained by the
H\"older inequality. For that purpose we need some control of
$\varphi$. From Lemma 2.3 it follows easily that for $w_+ \in H^3$,
$\nabla g (|w_+|^2) \in H^4$ and in particular $\nabla^k g(|w_+|^2)\in
L^{\infty}$ for $0 \leq k \leq 3$. Together with Lemma 3.4, this
provides an estimate of $\parallel \partial_t^j \nabla^k \varphi
\parallel_r$ for $j = 0, 1$, for $k = 1,2$ and $r = \infty$  and for $k
= 3$ and $r = 6$. With that information available, we apply the
H\"older inequality according to the following rules~: \par

(i) all the explicit powers of $x$ are attached to $w_+$. In addition,
whenever there appears a factor $\partial_t \varphi$ (with no space
derivative), one power $x$ is extracted from the $\widetilde{A}_1$ part
of $\partial_t \varphi$ and attached to $w_+$ (since $\widetilde{A}_1$
belongs to $L^{\infty}$ but a priori $\partial_t \varphi$ does not).
\par

(ii) The $x$ amputated contribution of $\partial_t \varphi$ generated
by rule (i) and all the factors $\partial_t^j \nabla^k \varphi$ with $k
= 1,2$ are estimated in $L^{\infty}$. The factors $\nabla^3\varphi$ are
estimated in $L^6$ (in fact in $H^1$). Such factors occur only from the
$t^{-1} x^s \nabla \Delta$ terms in the proof of (3.44).\par

(iii) The previous rules generate norms of the type $\parallel x^s
\nabla^k w_+ \parallel_r$ for $w_+$. Those norms are estimated by $H^1$
norms of the same quantities for $2 < r \leq 6$ and by $H^2$ norms for
$6 < r \leq \infty$.\par

(iv) The time dependence of the various terms follows from the explicit
$t$ dependence of the operators $Z$ of the previous list, together with
the fact that $\parallel \partial_t^j  \nabla^k \varphi  \parallel_r$
generates a factor $t^{-1}$ for $j = 1$ and a factor $\ell n\ t$ for $j
= 0$.\par

With the previous rules available, the proof reduces to an elementary
book keeping exercise, which will be omitted. We simply remark that the
dominant terms as regards $w_+$ have $x^3 \nabla$, $x^2 \nabla^2$ and
$x \nabla^3$ which are exactly controlled by the assumption $w_+ \in
H^{1,3} \cap H^{3,1}$, equivalent to the assumption $u_+ \in H^{3,1}
\cap H^{1,3}$. As regards the time dependence, the dominant terms come
from $x^s \nabla \varphi$ in the proof of (3.43) thereby generating a
factor $\ell n\ t$, and from $x^s \Delta \exp (-i \varphi )$ in the
proof of (3.44), generating $x^s |\nabla \varphi|^2$ and therefore a
factor $(\ell n\ t)^2$.\par

Finally (3.40) is the special case $j = k = 0$, $Z = 1$, $r = 3$ of
(3.45), while (3.41) follows from the estimate
$$\parallel \nabla u_a(t) \parallel_4 \ \leq \ \left ( \parallel xw_+
\parallel_4 \ + \ t^{-1} \left ( \parallel \nabla w_+ \parallel_4\ + \
\parallel \nabla \varphi \parallel_{\infty} \ \parallel w_+
\parallel_4\right ) \right ) t^{-3/4} \ . \eqno(3.46)$$ 
\par \nobreak \hfill $\sq$ \par

We can now complete the proof of Proposition 1.1.\\

\noi {\bf Proof of Proposition 1.1.} It suffices to show that the
assumptions of Proposition 2.2 are satisfied for the choice $h(t) =
t^{-1} (2 + \ell n\ t)^2$ made in Proposition 1.1. Now the assumptions
(\ref{2.17e}) (\ref{2.19e}) follow from (3.39) (3.42) of Proposition 3.1,
the assumption (\ref{2.20e}) follows from Lemma 3.3 and Corollary 3.1.
The assumption (\ref{2.21e}) follows from Lemma 3.2 as regards
$R_{1,2}$ and from (3.44) of Proposition 3.1 as regards $R_{1,1}$. The
assumptions (3.11) of Lemma 3.1 and (\ref{3.19e}) of Lemma 3.2 are
special cases of (3.43). The assumption (3.21) of Lemma 3.2 follows
from Lemma 3.3 and from the fact that if $A_0$ is a solution of the
free wave equation $\sq A_0$ in the Coulomb gauge $\nabla \cdot A_0 =
0$, with initial data $(A_+, \dot{A}_+)$, then also $x\cdot A_0$ is a
solution of the free wave equation, namely $\sq (x \cdot A_0) = 0$,
with initial data $(x \cdot A_+, x\cdot \dot{A}_+)$. Finally the
assumptions (\ref{2.22e}) (\ref{2.23e}) follow from Lemma 3.1. \par

The smallness conditions needed for Proposition 2.2 bear on $c_3$ and
$c_4$. Now from (3.40), $c_3 = \parallel w_+ \parallel_3$ while from
(3.41) or (3.46) 
$$c_4 =\ \parallel xw_+ \parallel_4 \ +  O\left ( T_0^{-1}
\ell n\ T_0 \right ) \ .$$

\noi Since the estimates are used only for $t \geq T$, one can replace
$T_0$ by $T$ in that expression, and the last term can be made
arbitrarily small by taking $T$ sufficiently large, so that the
smallness condition of $c_4$ reduces to the smallness of $\parallel
xw_+\parallel_4$.\par \nobreak \hfill $\sq$ \par

\noi {\bf Remark 3.1.} The regularity assumptions on $u_+$ or $w_+$
could be somewhat weakened. The strongest assumptions come from the
$\Delta w_+$ term in $R_{1,1}$ and from the $x^3$, $x \nabla^2$ and
$x^2\nabla$ operators $Z$ in the estimate of $\partial_t \nabla u_a$.
On the one hand the $\Delta w_+$ term in $R_{1,1}$ could be eliminated
by the choice 
$$w(t) = U(1/t)^* w_+$$

\noi at the expense of generating either a more complicated and less
explicit $\varphi$ or additional terms in $R_2$. On the other hand, we
have obtained on $L^6$ estimate of $\partial_t \nabla u_a$ whereas an
$L^4$ estimate was sufficient. Only a minor weakening of the
assumptions on $u_+$ could be achieved along those lines, and we shall
not press that point any further.


\newpage
\begin{thebibliography}{99}


\bibitem{1r} Ginibre, J., Ozawa, T. : Long range scattering for nonlinear Schr\"odinger and Hartree equations in space dimensions $n \geq 2$, Commun. Math. Phys. {\bf 151} (1993), 619-645.

\bibitem{2r} Ginibre, J., Velo, G. : Scattering theory in the energy space for a class of nonlinear Schr\"odinger equations, J. Math. Pures Appl. {\bf 64} (1985), 363-401.

\bibitem{3r} Ginibre, J., Velo, G. : Generalized Strichartz
inequalities for the wave equation, J. Funct. Anal. {\bf 133} (1995), 50-68.

\bibitem{4r} Ginibre, J., Velo, G. : Long range scattering and modified
wave operators for the Wave-Schr\"odinger system, Ann. H.P. {\bf 3} (2002),
537-612.

\bibitem{5r} Ginibre, J., Velo, G. : Long range scattering and modified
wave operators for the Wave-Schr\"odinger system II, Ann. H.P. {\bf 4} (2003),
973-999.

\bibitem{6r} Ginibre, J., Velo, G. : Long range scattering and modified
wave operators for the Maxwell-Schr\"odinger system I. The case of
vanishing asymptotic magnetic field, Commun. Math. Phys. {\bf 236} (2003),
395-448.

\bibitem{7r} Ginibre, J., Velo, G. : Scattering theory for the
Schr\"odinger equation in some external time dependent magnetic fields,
preprint, math.AP/0401355.

\bibitem{8r} Ginibre, J., Velo, G. : Long range scattering for the
Wave-Schr\"odinger system with large wave data and small Schr\"odinger data, preprint, math. AP/0406608.

\bibitem{9r} Guo, Y., Nakamitsu, K., Strauss, W. : Global finite energy solutions of the Maxwell-Schr\"odinger system, Commun. Math. Phys. {\bf 170} (1995), 181-196.

\bibitem{10r} Keel, M., Tao, T. : Endpoint Strichartz estimates, Amer. J. Math. {\bf 120} (1998), 955-980.

\bibitem{11r} Nakamitsu, K., Tsutsumi, M. : The Cauchy problem for the coupled Maxwell-Schr\"odinger equations, J. Math. Phys. {\bf 27} (1986), 211-216.

\bibitem{12r} Nakamura, M., Wada, T. : Local wellposedness for the Maxwell-Schr\"odinger equations, preprint, math. AP/0304486.

\bibitem{13r} Ozawa, T. : Long range scattering for nonlinear
Schr\"odinger equations in one space dimension, Commun. Math. Phys.
{\bf 139} (1991), 479-493.

\bibitem{14r} Ozawa, T., Tsutsumi, Y. : Asymptotic behaviour of solutions for the coupled
Klein-Gordon-Schr\"odinger equations, in Spectral and Scattering Theory and Applications, 
Adv. Stud. in Pure Math., Jap. Math. Soc. {\bf 23} (1994), 295-305.

\bibitem{15r} Shimomura, A. : Wave operators for the coupled
Klein-Gordon-Schr\"odinger equations in two space dimensions, 
Funkcial. Ekvac. {\bf 47} (2004), 63-82.

\bibitem{16r} Shimomura, A. : Scattering theory for the coupled
Klein-Gordon-Schr\"odinger equations in two space dimensions, J.
Math. Sci. Univ. Tokyo {\bf 10} (2003), 661-685.

\bibitem{17r} Shimomura, A. : Scattering theory for the coupled
Klein-Gordon-Schr\"odinger equations in two space dimensions II, 
Hokkaido Math. J., in press.

\bibitem{18r} Shimomura, A. : Modified wave operators for the coupled
Wave-Schr\"odinger equations in three space dimensions, Disc. Cont.
Dyn. Syst. {\bf 9} (2003), 1571-1586.

\bibitem{19r} Shimomura, A. : Modified wave operators for
Maxwell-Schr\"odinger equations in three space dimensions, Ann. H.P.
{\bf 4} (2003), 661-683.

\bibitem{20r} Strauss, W. : Non linear Wave equations, CMBS
Lecture Notes 73, Am. Math. Soc., Providence 1989.

\bibitem{21r} Tsutsumi, Y. : Global existence and asymptotic behaviour
for the Maxwell-Schr\"odinger system in three space dimensions, Commun.
Math. Phys. {\bf 151} (1993), 543-576.

\bibitem{22r} Yajima, K. : Existence of solutions for Schr\"odinger evolution equations,
Commun. Math. Phys. {\bf 110} (1987), 415-426.


\end{thebibliography}
\end{document}